\documentclass[11pt]{article}

\usepackage{amsmath,amsthm,amssymb}
\usepackage[english]{babel}
\usepackage{enumerate}
\usepackage[bookmarks, citecolor = blue]{hyperref}
\usepackage{algorithm,algpseudocode}

\newcommand{\Gr}{Gr\"obner }
\newcommand{\Pl}{Pl\"ucker }

\newcommand{\Proj}{\textnormal{Proj}\,}
\newcommand{\Spec}{\textnormal{Spec}\,}
\newcommand{\Supp}{\textnormal{Supp}\,}
\newcommand{\rk}{\mathnormal{rk}\,}

\newcommand{\Lm}{\textnormal{\footnotesize LM}}
\newcommand{\Lt}{\textnormal{\footnotesize LT}}
\newcommand{\ed}{\textnormal{ed}\,}
\newcommand{\V}{\mathcal{V}}
\newcommand{\Syz}{\textnormal{Syz}}
\newcommand{\sat}{{\textnormal{sat}}}

\newcommand{\id}{\mathfrak}
\newcommand{\idj}{\mathfrak{j}} 
\newcommand{\idh}{\mathfrak{h}} 
 
\newcommand{\St}{\mathcal{S}t}
\newcommand{\hilb}{{\mathcal{H}\textnormal{ilb}}}
\newcommand{\hilbp}{{\mathcal{H}\textnormal{ilb}_{p(z)}^n}}

\newcommand{\PP}{\mathbf{P}} 
\newcommand{\Af}{\mathbb{A}}
\newcommand{\ZZ}{\mathbb{Z}}
\newcommand{\T}{\mathbb{T}}
\newcommand{\HH}{\mathcal{H}}

 \numberwithin{equation}{section}
\newtheorem{lemma}{Lemma}[section] 
\newtheorem{theorem}[lemma]{Theorem} 
\newtheorem{cor}[lemma]{Corollary} 
\newtheorem{prop}[lemma]{Proposition} 
\newtheorem{cri}[lemma]{Criterion}

\newtheorem*{introThA}{Theorem A}
\newtheorem*{introThB}{Theorem B}
\newtheorem*{introThC}{Theorem C}

\theoremstyle{definition}
\newtheorem{notation}[lemma]{Notation} 
\newtheorem{remark}[lemma]{Remark} 
\newtheorem{defin}[lemma]{Definition} 
\newtheorem{example}[lemma]{Example} 

\DeclareMathAlphabet{\mathpzc}{OT1}{pzc}{m}{it} 

\begin{document} 

\title{Rational components of Hilbert schemes} 
\author{Paolo Lella\footnotemark[1] \and Margherita Roggero\footnote{Supported by PRIN (Geometria della variet\`a algebriche e dei loro spazi di moduli) co-financed by MIUR (2008).}} 
\date{\empty} 

\footnotetext{Mathematics Subject Classification 2000: 14C05. Keywords:  Hilbert schemes, Initial ideals}

\maketitle

\baselineskip=+.5cm

\begin{abstract} 
The Gr\"obner stratum of a monomial ideal $\id{j}$ is an affine  variety that parameterizes the  family of all ideals having $\id{j}$ as initial ideal (with respect to a fixed term ordering). The Gr\"obner strata can be equipped in a natural way with a structure of homogeneous variety and are in a close connection with Hilbert schemes of subschemes in the projective space $\PP^n$. Using properties of the Gr\"obner strata we prove some sufficient conditions for the rationality of components of $\hilbp$. We show for instance that all the smooth, irreducible components in $\hilbp$ (or in its support) and the Reeves and Stillman component $H_{RS}$ are rational. We also obtain  sufficient conditions  for isomorphisms between  strata corresponding to   pairs of ideals defining a same subscheme, that can strongly improve an explicit computation of their equations.
\end{abstract}

\section{Introduction}

The aim  of the present paper is to investigate effective methods to study the Hilbert scheme of subschemes in the projective space $\PP^n$, both on the theoretical and the computational point of view, using Gr\"obner basis tools. Several authors have been working in this direction during   last years (for instance \cite{carra,robbiano-2008}), but our motivations mainly refer to some ideas and hints contained in the paper \cite{NS} by Notari and Spreafico.
In order to get a stratification of $\hilbp$, they introduce some affine varieties $\St(\id{j})$ (here called Gr\"obner strata) parameterizing families of ideals in $k[X_0,\dots,X_n]$ having the same initial ideal $\id{j}$ with respect to a fixed term ordering on the monomials. When only homogeneous ideals in $k[X_0,\dots,X_n]$ are concerned, we write $\St_h(\idj)$. The ideal defining a Gr\"obner stratum springs out from a procedure based on Buchberger's algorithm, but involves a reduction with respect to a set of polynomials which is not a Gr\"obner basis. 

\medskip

It is not difficult to realize that the support of $\St(\id{j})$ only depends on the initial data (the term ordering, the ideal $\id{j}$, etc.), but one cannot be beforehand sure that  different choices in the reduction steps always lead to the same ideal. In other words it is  not clear if Gr\"obner strata are scheme-theoretically well defined. This is a crucial point that is underlined for instance in \cite[Section 3]{robbiano-2008}. In fact, if anyone wants to deduce properties of an Hilbert scheme using a Gr\"obner stratum, it is necessary to  consider carefully the non-reduced structures, because $\hilbp$ can have non-reduced components (see \cite{iarrobino1,Kleppe,Mum}).
 A first achievement in this paper is the following result (see Theorem \ref{th:cor-unica}):

\begin{introThA} The ideal defining $\St(\id{j})$ does not depend on the reduction choices. 
\end{introThA}

In fact we exhibit an equivalent, but intrinsic definition for the ideal of $\St(\id{j})$, which by the way also allows a great simplification in the procedure for an explicit   computation of this ideal.

\medskip 

A meaningful property that all Gr\"obner strata enjoy is that they are homogeneous with respect to some non-standard grading. Homogeneous varieties are of a very special type: for instance they can be isomorphically embedded in the Zariski tangent space at the origin, which is of course the smallest affine space in which such an embedding can be done. Therefore a smooth homogeneous variety has to be isomorphic to an affine space; in fact, the variety has the same dimension as the space in which it can be embedded.  Moreover one can obtain directly the ideal of $\St(\id{j})$ in the \lq\lq minimal embedding\rq\rq\ in the Zariski tangent space, through a preliminary detection of a maximal set of \lq\lq eliminable variables\rq\rq\ (we briefly resume this method in \S\ \ref{sec:hom_variety}). This is a second key point in our work, because one of the main difficulties usually met studying  Gr\"obner strata (and even more Hilbert schemes) is due to the huge number of variables that their equations involve, even in very simple cases. Besides the obvious computational gain, we would like to underline the interesting theoretical outcome of this method: most of our proofs are obtained just using the minimal embedding. 

A second useful tool that can simplify the computation of equations defining a stratum is given in Theorem \ref{th:saturato}. Though the strata corresponding to  ideals that define a same subscheme are in general non-isomorphic, however we show that  two of them give rise to isomorphic strata when they satisfy suitable sufficient condition, so that we can equivalently take into consideration the most convenient one. 

\begin{introThB} Let $\id{j}$ be  monomial ideal in $k[X_0, \dots, X_n]$ which is saturated and   Borel-fixed w.r.t. the order on the set of variables $X_0\succ X_1 \succ \dots \succ X_n$. 
\begin{enumerate}[i)]
\item For every $m$, there is a set of eliminable variables of the ideal 
defining $\St_h(\id{j}_{\geqslant m})$, that contains all variables except at most the ones 
appearing in polynomials $F$ such that $\Lt(F)=X^{\alpha}X_0^{m-\vert \alpha\vert}$,
where $X^\alpha$ is a minimal generator of $\id{j}$.
\item If  $X_{n-1}$ does not appear in any monomial of degree $m+1$ in the monomial basis of $\id{j}$, then $\St_h(\id{j}_{\geqslant m-1})\simeq \St_h(\id{j}_{\geqslant m})$. 
\item Especially, if $s$ is the maximal degree of monomials in the monomial basis of $\id{j}$ containing $X_{n-1}$, then $\St_h(\id{j}_{\geqslant s-1})\simeq \St_h(\id{j}_{\geqslant m})$ for every $m\geq s$.
\end{enumerate}
\end{introThB}

In \S\ \ref{sec:regularity} and \S\ \ref{sec:opensubset} we investigate more closely the natural connection between the Gr\"obner stratum $\St_h(\id{j})$ and the Hilbert scheme $\hilbp$, where $p(z)$ is the Hilbert polynomial of $k[X_0, \dots, X_n]/\id{j}$. As for every ideal $\id{i}$ whose initial ideal is $\id{j}$, the modules $k[X_0, \dots, X_n]/\id{i}$ and $k[X_0, \dots, X_n]/\id{j}$ share the same Hilbert function, there is an obvious set-theoretic inclusion $\St_h(\id{j}) \subseteq \hilbp$. However it is not a simple task to understand if this inclusion is an algebraic embedding or not. The paper \cite{NS} deals with the same question, but mainly concerns Gr\"obner strata of saturated ideals with respect to the term ordering $\mathtt{DegRevLex}$: note that every subscheme $Z$ in $\PP^n$ can be defined by the saturated ideal $I(Z)$. In this paper we prefer to consider a slightly different approach, modeled on the classical construction of the Hilbert schemes (see for instance \cite{B,HS}). For every $Z\in \hilbp$ we consider the ideal $I(Z)_{\geqslant r}$, where $r$ is the Gotzmann number of $p(z)$. As $r$ is the worst Castelnuovo-Mumford regularity for all $Z \in \hilbp$, the Gr\"obner strata (with respect to any term ordering) of   monomial ideals generated in degree $r$ cover $\hilbp$.

Moreover, the subset of $\hilbp$ corresponding to $\St_h(\id{j}_{\geqslant r})$ always contains the one corresponding to $\St_h(\id{j})$ and the inclusion can be strict, because points  in $\St_h(\id{j})$ correspond to ideals defining subschemes in $\PP^n$  with the same Hilbert \emph{function} as the subscheme $Z=\V(\id{j})$, while being in $\St_h(\id{j}_{\geqslant r})$ only requires the same Hilbert \emph{polynomial}. An interesting example of this type is that of the lexicographic saturated ideal $\id{L}$ such that $k[X_0, \dots, X_n]/\id{L}$ has Hilbert polynomial $p(z)$ and whose regularity is indeed the Gotzmann number $r$ of $p(z)$: in \S\ \ref{sec:lexsegment} we show that $\St_h(\id{L}_{\geqslant r})$ is isomorphic to an open subset of the Reeves and Stillman component $H_{RS}$ of $\hilbp$, while in general $\St_h(\id{L})$ corresponds to a locally closed  subscheme of   lower dimension (see \cite[Remark 4.8]{RT}).

The main reason of our setting is contained in Theorem \ref{th:aperto}. Let $p(z)$ be any admissible Hilbert polynomial in $\PP^n$ with Gotzmann number $r$ and let   $\prec$ be a fixed term ordering on monomials of $k[X_0, \dots, X_n]$. Following the classical construction, we consider  $\hilbp$ as a subscheme of a projective space through the Pl\"{u}cker embedding of the grassmannian $\mathbb{G}(t,M)$, where $M=\dim_k (k[X_0, \dots, X_n]_r)$ and $t=M-p(r)$. The simple remark that the Pl\"{u}cker coordinates correspond to   sets of $t$ distinct monomials of degree $r$ (that we can write in decreasing order with respect to $\prec$), allows us to get a lexicographic total order on them.
 If $\id{j}_0$ is a monomial ideal generated in degree $r$ such that $k[X_0,\dots,X_n]/\id{j}_0$ has Hilbert polynomial $p(z)$ with  Gotzmann number $r$, then we show that $\St_h(\id{j}_0)$
is the locally closed subscheme of $\hilbp$ given by the conditions that the Pl\"{u}cker coordinate corresponding to the monomial basis of  $\id{j}_0$ does not vanish and the bigger ones vanish.

As a consequence we are able to prove that every irreducible and reduced component of $\hilbp$ (or of its support) has an open subset which is a homogeneous affine variety (with respect to a non-standard grading). Especially, if $\id{j}$ is generated by the $t$ largest degree $r$ monomials (we call it a $(r,\prec)$-\emph{segment ideal}), then  $\St_h(\id{j})$ is  naturally isomorphic to an open subset of $\hilbp$. Therefore we can easily deduce a few interesting properties of rationality for the components of Hilbert schemes (see Theorem \ref{th:aperto} \emph{\ref{it:aperto_iii})}, Corollary  \ref{th:segment01}, Corollary \ref{th:segment}, Corollary \ref{th:lexlex}):

\begin{introThC} Let $H$ be an irreducible component of $\hilbp$.
\begin{itemize}
   \item If $H$ is smooth, then it is rational. The same holds for its support  $\Supp H$.
   \item If $H$ contains a smooth point which corresponds to a $(r,\prec)$-segment ideal (where $\prec $ is any term ordering), then $H$ is rational. The same holds for $\Supp H$.
   \item The Reeves and Stillman component $H_{RS}$ of $\hilbp$ is rational.
\end{itemize}
\end{introThC}

The last item can be obtained as a direct consequence of the previous one, because the lexicographic saturated ideal $\id{L}$ corresponds to a smooth point in $H_{RS}$, as proved by Reeves and Stillman in \cite{RS}, and  $\id{L}_{\geqslant r}$ is a $(r,\mathtt{Lex})$-segment ideal. However, we can also get a new proof of this fact, not applying the quoted result by Reeves and Stillman, but proving that $\St_h(\id{L}_{\geqslant r})$ is isomorphic to an affine space using our method based on the minimal embedding (see Theorem \ref{th:lexstratumopenrational}).

\medskip

In \S\ \ref{sec:example} we present a pseudo-code description of the procedures based on our results, that can be   implemented using one of the several softwares for symbolic computation. With such a procedure we are able to write equations for some Gr\"obner strata corresponding to the Hilbert scheme $\hilb_{4z}^3$. In this way we  find  a computational confirmation of the results obtained by Gotzmann in \cite{gotzmann-2008}, namely that   $\hilb_{4z}^3$ has  two components of dimensions 23 and 16 respectively and also some improvements. In fact, we also show  that  both components of that Hilbert scheme are rational (because each of them   contains a  smooth point corresponding to a segment ideal), they have transversal intersection (studied using a third segment ideal) and finally that the forth segment ideal allowed by the Hilbert polynomial $p(z)=4z$ is  singular point whose stratum has dimension 23 and embedding dimension 27. 

\medskip

The paper is organized as follows.  
\S\ \ref{sec:notation} contains some general notation.

In \S\ \ref{sec:stratum_ideal}, we take up the construction of Gr\"obner strata made in \cite{RT} and prove that they are well defined (Theorem \ref{th:cor-unica}).

In \S\ \ref{sec:hom_variety} we discuss the main properties of Gr\"obner strata as homogeneous varieties with respect to a non-standard grading and we obtain some useful criterion in order to know when Borel-fixed monomial ideals with the same saturation define the same Gr\"obner stratum (Theorem \ref{th:saturato}).

In \S\ \ref{sec:regularity}, we focus our attention on   ideals generated in degree $r$, where $r$ is the Gotzmann number of their Hilbert polynomials, and prove that their Gr\"obner strata can be defined by minors of suitable matrices (Proposition \ref{th:matriceA}).

\S\ \ref{sec:opensubset} represents the heart of the work. We show that there is a close connection between the above quoted matrices defining Gr\"obner strata and those appearing in the classical construction of Hilbert schemes and obtain as a consequence the main results of the paper about rational components (Theorem \ref{th:aperto} \emph{\ref{it:aperto_iii})}).

Finally, in \S\ \ref{sec:lexsegment}, we prove the rationality of the Reeves-Stillman component $H_{RS}$ of $\hilbp$ using our method, based on the minimal embedding and in \S\ \ref{sec:example} we apply this same method in order to perform some explicit computations about $\hilb_{4z}^3$.

\section{Notation} \label{sec:notation}

Throughout the paper, we will consider the following general notation. 
\begin{enumerate}
\item During the construction of Gr\"obner strata, we work on a field $k$ of any characteristic, whereas when we study the Hilbert scheme, we will suppose that $k$ is algebraically closed.

\item\label{not:general_monomial} $k[X_0,\dots,X_n]$ is the polynomial ring in the set of variables $X_0,\dots,X_n$ that we will often denote by the compact notation $X$, so that $k[X]:=k[X_0,\dots,X_n]$; we will denote by $X^\alpha$ the generic monomial in $k[X]$, where $\alpha$ represents a multi-index $(\alpha_0,\dots,\alpha_n)$, that is $X^\alpha := X_0^{\alpha_0} \cdots X_n^{\alpha_n}$. 
$\id{j}$ will be a monomial ideal in $k[X]$ with basis $\{X^{\gamma_1},\ldots,X^{\gamma_t}\}$ and  $\Syz(\id{j})$ its  $k[X]$-module  of  syzygies. 

$X^\alpha \mid X^\gamma$ means that $X^\alpha$ divides $X^\gamma$, that is there exists a monomial $X^\beta$ such that $X^\alpha \cdot X^\beta = X^\gamma$. If such monomial does not exist, we will write $X^\alpha \nmid X^\gamma$.

$\prec$ will be a fixed term ordering on the set $\T_X$ of monomials in $k[X]$ and we always assume that $X_0 \succ \dots \succ X_n$.
 As the term order $\prec$ is fixed, we often omit to indicate it. Given a monomial $X^\alpha$, we refer with $\min(X^\alpha)$ as the smallest variable dividing the monomial, that is $\min(X^\alpha) = \min \{X_i\ \text{ s.t. } X_i \mid X^\alpha\}$.
 We will denote also by $\prec$ its extension to the multiplicative group of Laurent monomials
 $\overline{\T}_{X}$ and the corresponding total ordering on $\ZZ^{n+1}$  given by $\alpha \prec \beta \Leftrightarrow X^{\alpha}\prec X^{\beta}$.

For every polynomial $F$ in $k[X]$ (or $k[X,C]$, $k[C]$), $\Lt(F)$ is its leading term with respect to the fixed term ordering; in the same way, if $\id{a}$ is an ideal, $\Lt(\id{a})$ is its initial ideal.

\item We will introduce a second set of variables $C_{i\alpha}$ that we will denote with $C$. So $k[X,C]$ will be the polynomial ring in the variables $X$ and $C$ and $\T_{X,C}$ the corresponding set of monomials. The term ordering on $\T_{X,C}$ will be induced by the term ordering on $\T_{X}$ and it will be an elimination ordering of the variables $X$ that will coincide with $\prec$ on $\T_{X}$: so we will denote by the same symbol $\prec $ also this term orderings on $\T_{X,C}$ and its restriction to $\T_{C}$.

\item \label{not:X_coeff} Let $G$ be any polynomial in $k[C,X]$. An \emph{$X$-monomial} of $G$ is a monomial of $\T_X$ that appears in $G$ considered as a polynomial in the variables $X$ with coefficients in the ring $k[C]$; the \emph{$X$-coefficients} of $G$ are the elements of $k[C]$ that are coefficients of an $X$-monomial. Note that the $X$-coefficients are  polynomials, but not necessary monomials.

\item Given any subscheme $Z $ in $\PP^n$, we will denote by $\Supp Z$ its support and by $I(Z)$ the saturated ideal in $k[X]$ that defines $Z$. Given any ideal $\id{a}$, we will denote by $\V(\id{a})$ the affine scheme $\Spec(k[X]/\id{a})$.

\item \label{not:hilb_scheme} $\hilbp$ will denote the Hilbert scheme parameterizing all subschemes $Z $ in $\PP^n$ with Hilbert polynomial $p(z)$.  $r$ will be the  Gotzmann number of $p(z)$, that is the worst Castelnuovo-Mumford regularity among subschemes parameterized by $\hilbp$. When we write that an ideal $\id{i} \subset k[X]$ belongs to $\hilbp$, we will mean that $\id{i} $ is generated in degree $r$ and that the Hilbert polynomial of $\Proj k[X]/\id{i}$ is $p(z)$. By abuse of notation we will say that any such ideal $\id{i}$ has Hilbert polynomial $p(z)$ referring to the Hilbert polynomial of the quotient, even if the real Hilbert polynomial of $\id{i}$ is $\binom{z+n}{n} - p(z)$.
\end{enumerate}

\section{The ideal of a Gr\"obner Stratum}\label{sec:stratum_ideal}

Now we introduce the Gr\"obner strata and prove some properties, generalizing definitions and results of the paper \cite{RT}. 

\begin{defin}\label{def:tail} 
The \emph{tail} of  $X^{\gamma}$  with respect to $\id{j}$ (and to the fixed term ordering $\prec$) is the set of monomials:
\begin{equation}
T_{\gamma}^\prec = \big\{ X^\alpha \in \T_X \ \big|\ X^\alpha \prec X^{\gamma},\  X^\alpha \notin \id{j} \big\}
\end{equation}
\end{defin}

Every ideal $\id{i}$ such that $\Lt(\id{i})=\id{j}=(X^{\gamma_1},\dots,X^{\gamma_t})$ has a reduced Gr\"obner basis of the type $\{f_1,\dots,f_t\}$ where:
\begin{equation}\label{eq:polibasepartenza1}
  f_i = X^{\gamma_i} + \sum_{X^{\alpha} \in T_{\gamma_i}^\prec} {c}_{i\alpha} X^{\alpha}
\end{equation}
and ${c}_{i\alpha }\in k$, ${c}_{i\alpha }=0$ except finitely many of them. It is very natural to parameterize the family of all the ideals $\id{i}$ by the coefficients ${c}_{i\alpha}$; in this way it corresponds to a subset of $k^{T^\prec}$, where $T^\prec=   T_{\gamma_1}^\prec\times\cdots\times T_{\gamma_t}^\prec$.

In many interesting cases, $T_{\gamma_i}^\prec$ are finite sets and so $k^{T^\prec}$ is an affine space: this happens for instance if $\id{j}$ is a zero-dimensional ideal or if $\prec$ is a suitable term ordering; in other cases, for instance when only homogeneous ideals are concerned, $T^\prec$ can be infinite, but we can restrict our interest to a suitable finite subset. The following definition extends and includes all the previous cases.

\begin{defin} \label{def:codaridotta}  
Let us fix $T=\{T_1, \dots, T_t\}$ where $T_i$ is a finite subset of the tail of $X^{\gamma_i}$ with respect to $\id{j}$. 
We will denote by $\St(\id{j},T)$ the family of all ideals $\id{i}$ in $k[X]$ such that $\Lt(\id{i})=\id{j}$ and whose reduced Gr\"obner basis ${f}_1,\dots,{f}_t$ is of the type:
\begin{equation}\label{eq:polibasepartenza2}
  {f}_i = X^{\gamma_i} + \sum_{X^{\alpha} \in T_{i}} {c}_{i\alpha} X^{\alpha}.
\end{equation}
Moreover we will use the following special notation:
\begin{enumerate}[i)]
  \item\label{it:codaridotta_i} $\St(\id{j})$, if  $T_i=T_{\gamma_i}^\prec$ (of course only if  $T_{\gamma_i}^\prec$ are finite sets):   $\St(\id{j})$ parameterizes all the ideals $\id i$ such that $\Lt(\id{i})=\id{j}$. 
  \item\label{it:codaridotta_ii} $\St_h(\id{j})$, if $T_i$ is the subset of $T_{\gamma_i}^\prec$ of the monomials with the same degree as $X^{\gamma_i}$: $\St_h(\id{j})$ parameterizes all the homogeneous ideals $\id i$ such that $\Lt(\id{i})=\id{j}$.
\end{enumerate}
\end{defin} 

\begin{remark}\label{rk:termordersignificance}
It will be clear later that the term ordering  affects the construction of a \Gr stratum only because it states which monomials can belong to  the  tails; in fact   two different term orderings giving the same tails will lead to the same \Gr strata.
\end{remark}

Every ideal $\id{i}$ in the family $\St(\id{j},T)$  is uniquely determined by a point in the affine space $\Af^N$ ($N=\sum_i \vert T_i \vert$) where we fix coordinates  ${C}_{i\alpha}$ corresponding to the coefficients $c_{i \alpha}$ that appear in  (\ref{eq:polibasepartenza2}). The subset of $\Af^N$ corresponding to $\St(\id{j},T)$ turns out to be a closed algebraic set. More precisely, we will see how it can be endowed in a very natural way with a structure of  affine subscheme, possibly reducible or non reduced, that is we will see that it can be obtained as the subscheme of $\Af^N$ defined by an ideal $\id{h}(\id{j},T)$ in $k[C]$, where $C$ is the set of variables $C_{i\alpha}$.

\medskip

In the following, we refer to the terminology introduced in Notation \ref{not:X_coeff} for what concerns   the polynomials in $k[X,C]$.

\begin{defin} \label{def:procedura} 
We will denote by $\id{h}(\id{j},T)$ and  $L(\id{j},T)$ respectively any ideal in $k[C]$ that can be obtained in the following way. 
\begin{itemize}
  \item Let $\mathcal{B} = \{F_1,\ldots,F_t\}$ be the set of polynomials in $k[X,C]$ given by:
\begin{equation}\label{eq:polibasegenerica}
  F_i = X^{\gamma_i} + \sum_{X^{\alpha} \in T_i}  {C}_{i\alpha} X^\alpha.  
\end{equation}
  \item Consider any term order in $k[X,C]$ which is an elimination order for the variables $X$ and that coincides with $\prec$ for monomials in $\T_X$; there will  be no confusion if we denote it by the same symbol $\prec$. With respect to such a term order, the leading term of $F_i$ is $X^{\gamma_i}$.
  \item Fix the subset $P$ of $\{ (i,j) \ |\ 1\leqslant i < j \leqslant m \}$ corresponding to any set of generators for $\Syz(\id{j})$;
  \item  For every  $(i,j)\in P$, let $R_{ij}$ be a complete reduction of the $S$-polynomial $S(F_i,F_j)$ with respect to $\mathcal{B}$.
  \item  For every  $(i,j)\in P$, let $M_{ij}$ be a complete reduction of $S(F_i,F_j)$ with respect to $\id{j}$.
  \item $\id{h}(\id{j},T)$ is the ideal in $k[C]$ generated by the $X$-coefficients of the polynomials $R_{ij}$, $(i,j)\in P$.
  \item $L(\id{j},T)$ is the $k$-vector space in $\langle C\rangle$ generated by the $X$-coefficients of $M_{ij}$, $(i,j)\in P$.
\end{itemize}
\end{defin}

It is almost evident, that the definition of $\id{h}(\id{j},T)$ is nothing else than   Buchberger's characterization of Gr\"obner basis if we think to the $C_{i\alpha}$'s as constant in $k$ instead of variables. In fact the variables $C$ do not appear in the leading terms of $F_i$ and so their specialization in $k$ commutes with reduction with respect to $\mathcal{B}$. Thus $(\dots,{c}_{i\alpha},\dots)$ is a closed point in the support of   $\V (\id{h}(\id{j},T))$ in $\Af^N$ if and only if it corresponds to polynomials $f_1,\dots,f_m$ in $k[X]$ that are a Gr\"obner basis. 
Then the support of $\V (\id{h}(\id{j},T))$ is uniquely defined; however a priori the ideal $\id{h}(\id{i},T)$ could depend on the choices we perform computing it, that is on the choice  of the set $P$ of generators for $\Syz(\id j)$ and on the choice of a reduction for the $S$-polynomials $S(F_i,F_j)$ with respect to $\mathcal{B}$ (which in general is not uniquely determined).  

Thanks again to Buchberger's criterion, we can prove that in fact $\id h(\id{j}, T)$ only depends on $\id{j}$, $T$ and of course $\prec$ because it can be defined in an equivalent intrinsic way. 

\begin{prop}\label{th:unica} Let  $\id{j} \subseteq k[X]$,  $\mathcal{B}=\{F_1,\dots, F_t\}\subset k[X,C]$ and $\prec$ be as above and consider an ideal $\id{a}$ in $k[C]$  with  Gr\"obner basis $\mathcal{A} $. The following are equivalent: 
\begin{enumerate}[i)]
\item\label{it:unica_1} $\mathcal{B} \cup \mathcal{A} $ is a  Gr\"obner basis in $k[X,C]$;
\item\label{it:unica_1'} $\id{a}$ contains the $X$-coefficients of all the polynomials in the ideal  $(F_1,\dots, F_t) k[X,C]$ that are reduced modulo $\id{j}$;
\item\label{it:unica_2} $\id{a}$  contains  all the $X$-coefficients   of every complete  reduction of $S(F_i,F_j)$ with respect to $\mathcal{B}$ for every $i,j$;  
\item\label{it:unica_3} $\id{a}$ contains  all the $X$-coefficients  of  some  (even partial) reduction  with respect to $\mathcal{B}$ of $S(F_i,F_j)$ for every $i,j$;
\item\label{it:unica_4} $\id{a}$ contains  all the $X$-coefficients  of  some  (even partial) reduction  with respect to $\mathcal{B}$ of $S(F_i,F_j)$, for every $(i,j)$ corresponding to a set of generators  of $\Syz(\id{j})$.
\end{enumerate}
\begin{proof} $\emph{\ref{it:unica_1})} \Rightarrow \emph{\ref{it:unica_1'})}$: let $G$  a polynomial in $(F_1,\dots, F_t) k[X,C]$ which is reduced modulo $\id{j}$. By  hypothesis, $G$ must be reducible to 0 through $\mathcal{B} \cup \mathcal{A}$, so that the next step of reduction have to be performed just using $\mathcal{A}$. But any step of reduction through $\mathcal{A}$ does not change the $X$-monomials and only modifies the $X$-coefficients; then  $G \stackrel{\mathcal{A}}{\longrightarrow} 0$, that is every $X$-coefficient in $G$ can be reduced to $0$ using $\mathcal{A}$: this shows that all the $X$-coefficients in $G$ belong to  $\id{a}$.

$\emph{\ref{it:unica_1'})} \Rightarrow \emph{\ref{it:unica_2})}$, $\emph{\ref{it:unica_2})} \Rightarrow \emph{\ref{it:unica_3})}$ and $\emph{\ref{it:unica_3})} \Rightarrow \emph{\ref{it:unica_4})}$ are obvious.

$\emph{\ref{it:unica_4})} \Rightarrow \emph{\ref{it:unica_1})}$: we can check that  $\mathcal{B} \cup \mathcal{A} $ is a Gr\"obner basis using the refined Buchberger criterion (see for instance \cite[Theorem 9, pag. 104]{CLO}). If $\mathcal{A} =\{a_1,\dots,a_r\}$, a set of generators for  $\Syz(X^{\gamma_1}, \dots, X^{\gamma_t}, \Lt(a_1), \dots, \Lt(a_r))$ can be obtained as the union of a set of generators for  $\Syz(X^{\gamma_1}, \dots, X^{\gamma_t})$, a set of generators for $ \Syz(\Lt(a_1), \dots, \Lt(a_r))$ and the obvious syzygies of $(X^{\gamma_i}, \Lt(a_j))$. Then: 
\begin{itemize}
 \item $S(a_i,a_j) \stackrel{\mathcal{B} \cup \mathcal{A}}{\longrightarrow} 0$, since $\mathcal{A}$ is a Gr\"obner basis and $\mathcal{A} \subseteq \mathcal{B}\cup\mathcal{A}$;
 \item $S(a_i,F_j) \stackrel{\mathcal{B} \cup \mathcal{A}}{\longrightarrow}0$, since the leading terms of $a_i$ and $F_j$ are coprime and $a_{i},\ F_j \in \mathcal{B}\cup\mathcal{A}$;
 \item $S(F_i,F_j) \stackrel{\mathcal{B}\cup\mathcal{A}}{\longrightarrow} 0$ in at least one way, by hypothesis. \qedhere
\end{itemize} 
\end{proof}
\end{prop}

There are many ideals $\id{a}$ fulfilling the equivalent conditions of Proposition \ref{th:unica}: for instance we can consider the irrelevant maximal ideal in $k[C]$ or any ideal obtained accordingly with condition \ref{it:unica_3}. Moreover, if $\id{a}$ satisfies those conditions and $\id{a}' \supset \id{a}$, then also $\id{a}'$  does, and if the ideals $\id{a}_l$ satisfy the conditions, then also their intersection $\bigcap \id{a}_l$ does. As a consequence of these remarks we obtain the proof of the uniqueness of the ideal $\id{h}(\id{j},T)$ given by Definition \ref{def:procedura}.

\begin{theorem}\label{th:cor-unica}  Let $\id j$ and $T$ as above. Then:   
\begin{enumerate}[i)]
\item\label{it:cor-unica_i} $\id h(\id{j},T)$ is uniquely defined; in fact $\id h(\id{j},T)=\bigcap \id{a}$, $\id{a}$  satisfying the equivalent conditions of Proposition \ref{th:unica}
\item\label{it:cor-unica_ii} $L(\id{j},T)$ is uniquely defined.
\end{enumerate}
\begin{proof} \emph{\ref{it:cor-unica_i})}: $\id{h} $ is one of the ideals $\id{a}$, because it satisfies condition \emph{\ref{it:unica_4})}; on the other hand, if $\id{a}$ satisfies condition \emph{\ref{it:unica_2})}, then clearly $\id{a} \supseteq \id{h}$.

For \emph{\ref{it:cor-unica_ii})} it is sufficient to observe that the generators for $L(\id{i},T)$  are the degree 1 homogeneous components (here \lq\lq homogeneous\rq\rq\ is related to the usual grading of $k[C]$ that is the $\ZZ$-grading with variables of degree 1) of the generators of $\id{h}(\id{i},T)$ given in its construction (Definition \ref{def:procedura}).
\end{proof}
\end{theorem}

By abuse of notation we will denote by the same symbol $\St(\id{j},T)$  the family of ideals and the  subscheme in $\Af^N$ given by the ideal $\id{h}(\id{i},T)$.   Note that $\id{h}(\id{i},T)$ is not always a prime ideal and so $\St(\id{j},T)$  is not necessarily irreducible nor reduced, as the following trivial example shows.

\begin{example} Let $\idj=(x^2,xy)\subset k[x,y]$ and  $\prec$ be any term ordering.  Let us  choose $T=\big\{T_{x^2}=\emptyset,T_{xy}=\{y\}\big\}$ and  construct the ideal of the \Gr stratum $\St(\idj,T)$ according to Definition \ref{def:procedura}:  
\[
\{ F_1=x^2, \ F_2=xy+Cy\}, \ S_{12}=yF_1-xF_2=-Cxy\stackrel{\{F_1,F_2\}}\longrightarrow R_{12}=-Cxy+CF_2=C^2y.
\] 
Then $\idh(\idj,T)=(C^2)$ that is $\St(\idj,T)$ is a double point in the affine space $\mathbb{A}^1$.
\end{example}

\section{Gr\"obner strata are homogeneous varieties}\label{sec:hom_variety}

In this section we will see how every Gr\"obner stratum $\St(\id{j},T)$ is in a very natural way homogeneous with respect to a suitable non-standard grading on $k[C]$, so that we can apply the nice properties typical of this kind of schemes and especially those obtained in \cite{RT} and in \cite{FR}. 

For the meaning of $\id{j}$, $k[X]$, $\{X^{\gamma_1},\ldots,X^{\gamma_t}\}$ and $\prec$ we refer to Notation \ref{not:general_monomial} and for $k[X,C]$, $\{F_1, \dots, F_t\}$, $\St(\id{j},T)$, $\id{h}(\id{i},T)$ to the previous section.

First of all, we recall the definitions and properties that we will use more often.

\begin{defin}\label{def:lambda} We will consider  $k[X, C]$  and  $k[C]$ as graded ring over   the totally ordered group $(\ZZ^{n+1},+,\prec)$ with grading $\lambda$ given by  $\lambda(X^{\alpha})=\alpha$ and   $\lambda(C_{i\alpha})=\gamma_i -\alpha $. 
\end{defin}
 
 As we will use also the usual grading over $\ZZ$ where all the variables have degree 1, we will  always write explicitly   the symbol $\lambda$ when the above defined grading is concerned (so, $\lambda$-degree $l$ with $l\in \ZZ^{n+1}$, $\lambda$-homogeneous of degree $l$  etc.) and leave the simple terms when the usual grading is involved (so, degree $r$ with $r\in \ZZ$, homogeneous of degree $r$  etc.).

\begin{prop}\label{th:lambda1} (See \cite[Lemma 2.8]{RT})  
\begin{enumerate}[i)]
\item\label{it:lambda1_i} The grading $\lambda$ is positive. 
\item\label{it:lambda1_ii} $\id{h}(\id{j},T)$ is a $\lambda$-homogeneous ideal.
\end{enumerate}
\begin{proof} \emph{\ref{it:lambda1_i})}  Let us observe that all the variables have $\lambda$-degree higher than that of the constant 1. In fact  $\lambda(X_i) \succ \lambda(1)$ because $\prec$ is a term ordering and $\lambda(C_{i\alpha})\succ \lambda(1)$ because, $X^{\gamma_i}\succ X^{\alpha}$ by definition of tails. As well known, this condition is equivalent to the positivity of the grading (see \cite[Chapter 4]{KrRob2}).

\emph{\ref{it:lambda1_ii})} We observe that $\lambda$ on $\overline{\T}_{C}$ is the restriction of the grading on  $\overline{\T}_{X,C}$. Every monomial that appears in $F_i$ is of the type $C_{i\alpha}X^{\alpha}$ and so its $\lambda$-degree is $\lambda(C_{i\alpha}X^{\alpha})=\lambda(C_{i\alpha})+\lambda(X^{\alpha})=\gamma_i$.
Thus all the polynomials $F_i$ are $\lambda$-homogeneous and then also the $S$-polynomials $S(F_i,F_j)$ and their reductions are $\lambda$-homogeneous. Finally, the $X$-coefficients in any $\lambda$-homogeneous polynomial (which are polynomials in $k[C]$) are $\lambda$-homogeneous. 
\end{proof}
\end{prop}

We now recall some properties of $L(\id{j},T)$ (see also \cite[Proposition 2.4]{RT} and \cite[Theorem 3.2]{FR}).

\begin{prop}\label{th:richiami} 
The linear space $\V(L(\id{j},T))$ can be naturally identified with the Zariski tangent space to $\St(\id{j},T)$ at the origin.

If $C''\subset C$ is any subset of $\ed := \dim \V(L(\id{j},T))$ variables such that   $L(\id{j},T)\oplus\langle C'' \rangle=\langle C\rangle$, then $\id{h}(\id{j},T)\cap k[C'']$ defines a $\lambda$-homogeneous subvariety in $\Af^{\ed }$ isomorphic to $\St(\id{h},T)$.
\end{prop}

We may summarize the previous result saying that \emph{$\St(\id{h},T)$ can be embedded in its Zariski tangent space at the origin}. This explains the following terminology.

\begin{defin}\label{def:eliminabili} The number $\ed   $ is the \emph{embedding dimension} of $\St(\id{j},T)$. The complement  $C':=C\setminus C''$ is a \emph{maximal set of eliminable variables} for $\id{h}(\id{j},T)$.
\end{defin}

\begin{cor}\label{th:liscio} In the above notation, the following statements are equivalent:
\begin{enumerate}
 \item $\St (\id{j},T)\simeq \Af^{\ed }$;
  \item $\St(\id{j},T)$  is   smooth;
 \item the origin is a smooth point for $\St(\id{j},T)$;
 \item $\ed  \leqslant \dim \St( \id{j},T)$.
\end{enumerate}
\end{cor}

Note that in general a maximal set of eliminable variables (and so its complementary) is not uniquely determined. However, if $C_{i\alpha}\in L(\id{j},T)$, then $C_{i\alpha}$  belongs to any set of eliminable variables; on the other hand, if $C_{i\alpha}$ does not appear in any element of $L(\id{j},T)$, then $C_{i\alpha}$ does not belong to any set of eliminable variables.

\medskip

There is an easy criterion that allows us to decide if a variable is eliminable or not.
\begin{cri}\label{criterio} Let $\Lt(F_i)=X^{\gamma_i}$, $\Lt(F_j)=X^{\gamma_j}$ and let $C_{i\beta}$ be a variable appearing in the tail of $F_i$. Using the reduction with respect to $\id{j}$ of a $\lambda$-homogeneous polynomial $X^{\delta}F_i-X^{\eta}F_j$ we can see that:
\begin{enumerate}[i)]
\item\label{it:criterio_i} if $X^{\delta+\beta}\notin \id{j}$ and  $X^{\delta +\beta -\eta}$ is not a monomial that appears in $F_j$, then $C_{i\beta}\in L(\id{j},T)$;
\item\label{it:criterio_ii} if $X^{\delta+\beta}\notin \id{j}$ and $X^{\beta'}=X^{\delta +\beta -\eta}$ is a monomial that appears in $F_j$, then $C_{i\beta}-C_{j\beta'}\in L(\id{j},T)$
\end{enumerate}
Moreover if $C_{i\beta}-C_{j\beta'} \in L(\id{j},T)$, then every maximal set of eliminable variables must contain at least either one of them.
\end{cri}
In most cases the number $N = |C|$ is very big and $\id{h}(\id{j},T)$ needs a lot of generators so that finding it explicitly is a very heavy computation. On the contrary $L(\id{j},T)$ is very fast to compute and so we can easily obtain a set of eliminable variables $C'$; a forgoing knowledge of $C'$ allows a simpler computation of the ideal  $\id{h}(\id{j},T)\cap k[C\setminus C']$ that gives $\St(\id{j},T)$ embedded in the affine space of minimal dimension $\Af^{\ed }$.

Furthermore, in many interesting cases we can greatly bring down the number of involved variables thanks to another kind of argument.

\begin{theorem}\label{th:saturato} Let $\id{j} \subset k[X_0,\ldots,X_n]$ be a Borel-fixed saturated monomial ideal with  basis $B$,   $m$ any  integer and   $\id{h}_m :=\id{h}(\id{j}_{\geqslant m})$  the ideal  of $\St_h(\id{j}_{\geqslant m})$  as in Definition \ref{def:procedura}. 
\begin{enumerate}[i)]
\item\label{it:saturato_i}  There is a set of eliminable variables for $\id{h}_m$ 
that contains all variables except at most the ones 
appearing in polynomials $F_i$ whose leading term is   either  $X^\gamma \in B_{\geq m}$   or  $X^\alpha X_n^{m-\vert \alpha \vert } $, where $X^\alpha \in B_{<m}$.
\item\label{it:saturato_ii}  $\St_h(\id{j}_{\geqslant m-1})$ is a closed subscheme of  $  \St_h(\id{j}_{\geqslant m})$. More precisely   $\St_h(\id{j}_{\geqslant m-1})\simeq \St_h(\id{j}_{\geqslant m},T)$  where $T$ contains the complete tail of a monomial in the basis of $\id{j}_{\geqslant m}$ if it is not divided by $X_n$, and a tail containing only  monomials divided by $X_n$ otherwise.
\item\label{it:saturato_iii}  If  $X_{n-1}$ does not appear in any monomial of degree  $m+1$ in the monomial basis of $\id{j}$, then $\St_h(\id{j}_{\geqslant m-1})\simeq \St_h(\id{j}_{\geqslant m})$. 
\item\label{it:saturato_v}  If  $X_{n-1}$ appears in $N$ monomials of degree $m+1$ in the monomial basis of $\id{j}$, then $\ed \St_h(\id{j}_{\geqslant m})\geqslant \ed \St_h(\id{j}_{\geqslant m-1}) + NM$, where $M$ is the number of monomials of the basis of $\id{j}$ of degree smaller than $m+1$. 
\item\label{it:saturato_vv} $\St_h(\id{j}_{\geqslant m-1})\not \simeq \St_h(\id{j}_{\geqslant m})$ if and only if   $X_{n-1}$  appears in monomials of degree $m+1$ in the monomial basis of $\id{j}$ and  $\id{j}_{\geqslant m-1}\neq  \id{j}_{\geqslant m}$.
\item\label{it:saturato_iv} If $s$ is the maximal degree of a monomial divided by $X_{n-1}$ in the monomial basis of $\id{j}$, then $\St_h(\id{j}_{\geqslant s-1})\simeq \St_h(\id{j}_{\geqslant m})$ for every $m\geqslant s$.
\end{enumerate}

\begin{proof} \emph{\ref{it:saturato_i})} Let us consider any monomial $X^\eta$ in the monomial basis of $\id{j}_{\geqslant m}$ which does not belong to $B_{\geqslant m}$ and such it that could be written as $X^\eta=X^\alpha X^\epsilon$ where $X^\alpha$ is a minimal generator of $\id{j}$ of degree $d<m$ and $X^\epsilon$ is a monomial of degree $m-d$,  $X^\epsilon \neq X_n^{m-d}$. Then among the polynomials $F_i$ there are:
\[
 \begin{split}
&  F = X^\alpha X_n^{m-d} + \sum C_{\beta} X^\beta,\\
&  F' = X^{\alpha+\varepsilon} + \sum C'_{\delta} X^\delta. 
 \end{split}
\]
We have to prove that all the variables $C'$ that appear in $F'$ can be eliminated.
 
The $S$-polynomial of $F$ and $F'$ is:
\[
 S(F,F') = X_n^p F' - X^{\varepsilon'} F = 
\sum C'_{\delta} X^\delta X_n^p - \sum C_{\beta} X^{\beta+\varepsilon'}.
\]
No monomial $X^\delta X_n^p$ in the first summand belongs to $\id{j}_{\geqslant m}$ because   $X^\delta \notin \id{j}$  and $\id{j}$ is saturated and Borel-fixed. Thus, the linear part of the coefficient of $X^\delta X_n^p$ in  the reduction of this $S$-polynomial will be either  $C'_{\delta}$ or $C'_{\delta}-C_{\beta} $.  Then $C'$ is a set of eliminable   variables  for $\id{j}_{\geqslant m}$.

\medskip

\emph{\ref{it:saturato_ii})} The first part of this statement is  a special case of general facts proved in \cite[\S 3]{HS}.

We directly prove the  second part  (which  implies the first one). Here  we  denote by $X^\alpha$ and  $X^\gamma$   the monomials in the basis of $\id{j}_{\geqslant m-1}$ of degree $m-1$ and $\geq m$ respectively, and we set:
\[\begin{split} G_{\alpha} & := X^\alpha  + \sum C_{\alpha\delta} X^\delta   \\
  G_{\gamma} & := X^\gamma  + \sum C_{\gamma\eta} X^\eta   
  \end{split} \]
  where  $ X^\delta $ varies among all   monomials of degree $m-1$ in the tail of $ X^\alpha$ and $X^\eta $ among those of the same degree as $ X^\gamma$  in its tail.   Applying the  procedure described in Definition \ref{def:procedura} on the set of polynomials $G$ we  define $\St_h(\id{j}_{\geqslant m-1})$   by an ideal $\overline{\id{h}} \subset k[C]$.
  
  The basis of $\id{j}_{\geqslant m}$ is made by  monomials of the following three types:
\begin{itemize}
\item monomials $X^\gamma$ of degree $\geq m$, that  also belong to   the basis of $\id{j}_{\geqslant m-1}$;
\item   monomials $X^\alpha X_n$    such that  $X^\alpha $ is any monomial of degree $m-1$ in the basis of $\id{j}_{\geqslant m-1}$;
\item  monomials $X^\alpha X_i$    of degree $m$ such that $X^\alpha $ is as above and $ \min(X^\alpha)\geq X_i \neq X_n$.
\end{itemize} 
We set:
 \[ \label{listaF}
 \begin{split}
 F_{\alpha n} & := X^\alpha X_n + \sum C_{\alpha\delta} X^\delta X_n \\
 F_{\alpha i} & := X^\alpha X_i + \sum C'_{\alpha i \tau} X^\tau  \quad \vert \tau\vert= m    \quad X^\tau \prec X^\alpha X_i \\
 F_{\gamma} & := X^\gamma  + \sum C_{\gamma\eta} X^\eta
 \end{split}\]
  Note that we use  the same  names   for some of the coefficients that appears in  polynomials $F$ and $G$, so that $F_{\alpha n}=X_n G_\alpha$ and $F_{\gamma}=G_{\gamma}$.
 Applying the  procedure described in Definition \ref{def:procedura} on the set of polynomials $F$  we obtain an ideal  $\id{h}' \subset k[C,C']$ defining $\St_h(\id{j}_{\geqslant m},T)$. 
  
    Thanks  to i) we know that $C'$ is a set of eliminable variables for $\id{h}'$ and so $\St_h(\id{j}_{\geqslant m},T)$ is also defined by $\id{h}=\id{h}'\cap k[C]$. The statement follows once we show that $\overline{\id{h}} =\id{h}$.
  
\medskip 
In order to eliminate the variables $C'$ we consider every monomial   $X^\alpha X_i =\Lm (F_{\alpha i} )$ and reduce it   using the polynomials $G$. In this way we obtain  a polynomial $H_{ \alpha i} \in (G)k[X,C]$ such that  $X^\alpha X_i+H_{ \alpha i} $ is completely reduced w.r.t. $\id{j}$.  Then    also  $X^\alpha X_iX_n+H_{ \alpha i}X_n +\sum C'_{\alpha i \tau} X^\tau X_n   $ (i.e. $F_{\alpha i}X_n+H_{ \alpha i}X_n  $) is  reduced modulo $\id{j}$ and moreover it belongs to $(F)k[X,C,C']$  because $X_nG\subseteq (F)k[X,C,C']$.    Its $X$-coefficients belong to $\id{h}'$, because the ideal $\id{h}'$  is generated by the  $X$-coefficient of  the polynomials in $   (F)k[X,C,C']$) that are reduced modulo $\id{j}_{m-1}$ or modulo $\id{j}$, which is the same (Proposition \ref{th:unica} \ref{it:unica_1'} and Theorem \ref{th:cor-unica}). The $X$-coefficients of  $F_{\alpha i}X_n+H_{ \alpha i}X_n  $ are also the $X$-coefficients of $F_{\alpha i}+H_{ \alpha i} $, and are precisely the set of polynomials of the type $C'_{ \alpha i \tau}-\phi_{ \alpha i \tau}(C)$ that allow us to   eliminate   the variables $C'$. So  the elimination of   $C'$ is obtained simply  putting $C'_{ \alpha i \tau}=\phi_{ \alpha i \tau}(C)$. In this way $F_{\alpha i}$ becomes $-H_{ \alpha i}$ that belongs to $(G)k[X,C]$.

  The ideal $\id{h}$, obtained from $\id{h}'$ eliminating $C'$, can  also be obtained first eliminating $C'$ and after taking $X$-coefficients, because the procedure of eliminating   $C'$ and that of taking $X$-coefficients commute.  So   $\id{h}$ is generated by the $X$-coefficients of polynomials in  $(X_nG_\alpha,-H_{ \alpha i}, G_\gamma)k[X,C]$ that are reduced modulo $\id{j}$. 
  
  Hence $\id{h}\subseteq \overline{\id{h}}$ because $(X_nG_\alpha,-H_{ \alpha i}, G_\gamma)k[X,C]\subset (G)k[X,C]$.

 On the other hand, $X_n(G)k[X,C]=(X_nG_\alpha, X_nG_\gamma)k[X,C]\subset (X_nG_\alpha, G_\gamma)k[X,C]$. Moreover two polynomials $Q$ and $X_nQ$ have the same  $X$-coefficients  and either one is   reduced modulo $\id{j}$ if and only the other is. Hence  we obtain the opposite inclusion  $\overline{\id{h}}\subseteq \id{h}$ and  conclude.

  \medskip
  
\emph{\ref{it:saturato_iii})} We use \emph{\ref{it:saturato_ii})} and prove that in the present hypothesis, $\St_h(\id{j}_{\geqslant m})\simeq \St_h(\id{j}_{\geqslant m},T)$, where $T$ is defined as in  \emph{\ref{it:saturato_ii})}.    Following Definition \ref{def:procedura}, we obtain the ideal  $\id{h}_m$ of $\St_h(\id{j}_{\geqslant m})$  using:
\[
 \begin{split}
 F''_{\alpha n} & := X^\alpha X_n + \sum C_{\alpha\delta} X^\delta X_n +\sum C''_{\alpha\sigma} X^\sigma \quad, \quad X_n \nmid X^\sigma\\
 F_{\alpha i} & := X^\alpha X_i + \sum C'_{\alpha i \tau} X^\tau   \\
 F_{\gamma} & := X^\gamma  + \sum C_{\gamma\eta} X^\eta= G_\gamma.
 \end{split}\]
 Note that $ F_{\alpha i}$ and $ F_{\gamma}$ are as in  \emph{\ref{it:saturato_ii})}, but   all the degree $m$ monomials of  the  tail of $X_nX^\alpha$ appear in $F''_{\alpha n}$, and not only those  divided by $X_n$.
  
For every monomial $X^\alpha$ of degree $m-1$ in the basis of $\id{j}_{m-1}$, let us consider  the $S$-polynomial:
$$S(F''_{\alpha n} ,F_{\alpha n-1})= \sum C_{\alpha\delta} X^\delta X_{n-1}X_n +\sum C''_{\alpha\sigma} X^\sigma X_{n-1} - \sum C'_{\alpha i \tau} X^\tau X_n.$$ 
 By our hypothesis no monomial appearing in it belongs to $\id{j}_{m}$. In fact  $X^\sigma X_{n-1} \in \id{j}$ if and only if it is a minimal generator of $\id{j}$, which is excluded by hypothesis because its degree is $m+1$, or it is of the type $X^\alpha X_a$ with $X^\alpha$ minimal generator of   $ \id{j}_m$  and  $X_a=\min (X^\sigma X_{n-1})=X_{ n-1}$, while $X^\sigma \notin \id{j}_m$. Then $S(F''_{\alpha n} ,F_{\alpha n-1})$ is already reduced with respect to $\id{j}_m$ and so  its $X$-coefficients  belong to $\id{h}_m$. Especially, as both $X^\delta X_{n-1}X_n$ and   $X^\tau X_n$ are multiple of $X_n$, while $X^\sigma X_{n-1}$ is not, the coefficient of $X^\sigma X_{n-1}$ is simply  $C''_{\alpha\sigma} $ so that each $C''_{\alpha\sigma} $ belongs to $\id{h}_m$.  Hence we can eliminate all the variables $C''$, just putting them equal to 0. In this way    $F''_{\alpha n}$ becomes $F_{\alpha n}$ as in (\ref{listaF}) and $\St_h(\id{j}_{\geqslant m})\simeq \St_h(\id{j}_{\geqslant m},T)$, where $T$ is as in \emph{\ref{it:saturato_ii})}, and we conclude  because $\St_h(\id{j}_{\geqslant m},T) \simeq \St_h(\id{j}_{\geqslant m-1}) $. 
 
\emph{\ref{it:saturato_v})} By \emph{\ref{it:saturato_ii})}, we know that $\ed \St_h(\id{j}_{\geqslant m}) \geqslant \ed \St_h(\id{j}_{\geqslant m},T) = \St_h(\id{j}_{\geqslant m-1})$, where the tails defined in $T$ contain only monomials divided by $X_n$. Let us now consider a monomial $X^\alpha$ among the generators of $\id{j}$ of degree smaller than $m+1$ and a generator $X^\gamma$ of degree $m+1$ divided by $X_{n-1}$. Computing the stratum $\St_h(\id{j}_{\geqslant m})$, in the tail of $X^\alpha X_n^{m-\vert\alpha\vert}$ there is the monomial $X^\beta = X^\gamma/X_{n-1}$ not belonging to $T$. Let us call $D$ the coefficient of $X^\beta$, that is
\[
F = X^\alpha X_n^{m-\vert\alpha\vert} + \ldots	 + D X^\beta + \ldots.
\]
Thinking about the syzygies of the ideal $\id{j}$, it is easy to see that in any $S$-polynomial, $F$ is surely multiplied by a monomial $X^\delta$ divided at least by one variable $X_i,\ i < n$. Therefore in every $S$-polynomial the monomial $X^\beta X^\delta = (X^\beta X_i) X^{\delta'}$ belongs to $\id{j}$ because of the Borel-fixed hypothesis, so that it can be reduced. Finally there is no equation involving the variable $D$, so it is free and it cannot be eliminated. Repeating the reasoning for the $M$ minimal generators of degree smaller than $m+1$ and for the $N$ generators divided by $X_{n-1}$ of degree $m+1$, we obtain the thesis.
  
\emph{\ref{it:saturato_vv})} straightforward applying \emph{\ref{it:saturato_v})}.  \emph{\ref{it:saturato_iv})} straightforward applying \emph{\ref{it:saturato_iii})}.
\end{proof}
\end{theorem}

With the following examples, we want to underline again the not so crucial role played by term ordering in this construction (Example \ref{ex:strato1}) and we want to show (Example \ref{ex:strato2} and Example \ref{rk:stratoLex}) that the estimate of growth of the embedding dimension of the stratum introduced in Theorem \ref{th:saturato} \emph{\ref{it:saturato_v})} is a lower bound.

\begin{example}\label{ex:strato1}
Let us consider the ideals $\id{i} = (X_0,X_1^2,X_1 X_2)$ and $\id{j} = \id{i}_{\geqslant 2} = (X_0^2,X_0 X_1,$ $X_0 X_2, X_0 X_3,X_1^2 , X_1 X_2)$ in the ring $k[X_0,X_1,X_2,X_3]$ and the strata of the ideal $\id{j}$ according to two different term orderings: $\St_h(\id{j},\mathtt{Lex})$ and $\St_h(\id{j},\mathtt{DegRevLex})$. In the first case there are at first $24$ new variables $C$, whereas in the second case they are $23$, so we may guess that the family of the ideals with initial ideal $\id{j}$ w.r.t. $\mathtt{Lex}$ could be different from the family of the ideals with initial ideal $\id{j}$ w.r.t. $\mathtt{DegRevLex}$.

However\hfill applying\hfill Theorem\hfill \ref{th:saturato},\hfill we\hfill can\hfill see\hfill that\hfill $\St_h(\id{j},\mathtt{Lex}) \simeq \St_h(\id{i},\mathtt{Lex})$\hfill and\\ $\St_h(\id{j},\mathtt{DegRevLex}) \simeq \St_h(\id{i},\mathtt{DegRevLex})$. Now the tails of the 3 monomials that generate $\id{i}$ are the same w.r.t. both term orders and then (see Remark \ref{rk:termordersignificance})
\[
 \St_h(\id{j},\mathtt{Lex}) \simeq \St_h(\id{i},\mathtt{Lex}) = \St_h(\id{i},\mathtt{DegRevLex}) \simeq \St_h(\id{j},\mathtt{DegRevLex}).
\]
\end{example}

\begin{example}\label{ex:strato2}
Let us consider the polynomial ring $k[X_0,X_1,X_2,X_3]$, the ideal $\id{j} =(X_0^2,$ $X_0 X_1,X_0 X_2^4,X_1^7,X_1^6 X_2^2)$ and any term ordering given by a matrix with first row  $(23,5,2,1)$. By the previous theorem we know that
\[
\St_h(\id{j}) \simeq \St_h(\id{j}_{\geqslant 3}),\quad \St_h(\id{j}_{\geqslant 4}) \simeq \St_h(\id{j}_{\geqslant 5}) \simeq \St_h(\id{j}_{\geqslant 6}), \quad \St_h(\id{j}_{\geqslant 7}) \simeq \St_h(\id{j}_{\geqslant m}),\ \forall\ m \geqslant 8
\]
and
\[
\ed \St_h(\id{j}_{\geqslant 4}) \geqslant \ed \St_h(\id{j}) + 2 \qquad \ed \St_h(\id{j}_{\geqslant 7}) \geqslant \ed \St_h(\id{j}_{\geqslant 4}) + 3
\]
By a direct computation, we find $\ed \St_h(\id{j}) = 46$,  $\ed \St_h(\id{j}_{\geqslant 4}) = 50$ and  $\ed \St_h(\id{j}_{\geqslant 7}) = 56$.
\end{example}

\begin{example}\label{rk:stratoLex}
There\hfill are\hfill at\hfill most\hfill two\hfill possible\hfill classes\hfill of\hfill isomorphism\hfill for\hfill the\hfill strata\\ $\St_h(\id{L}_{\geqslant m})$, where $\id{L}$ is a lexicographic ideal: $\St_h(\id{L})$ and $\St_h(\id{L}_{\geqslant r-1})$, where $r$ is
the maximal degree of a minimal generator, in fact the variable $X_{n-1}$ appears (if it does) only in the generator of degree $r$. Called $b$ the number of generators of degree $r$, applying Theorem \ref{th:saturato} \emph{\ref{it:saturato_v})}, we have
\[
\ed \St_h(\id{L}_{\geqslant r-1}) \geqslant \ed \St_h(\id{L}) + n-b.
\]

If the monomial of maximal degree in the basis does not contain the variable $X_{n-1}$, we have $\St_h(\id{L}_{\geqslant m}) \simeq \St_h(\id{L}),\ \forall\ m$. 
\end{example}

We conclude this section with a result similar to the one stated in Theorem \ref{th:saturato} that concerns only the case of homogeneous Gr\"obner strata w.r.t. $\mathtt{DegRevLex}$.
\begin{prop}
Let $\id{j}$ be a Borel-fixed saturated monomial ideal and let $\prec$ be the $\mathtt{DegRevLex}$ term ordering. Then
\begin{equation*}
 \St_h(\id{j}) \simeq \St_h(\id{j}_{\geqslant m}),\quad \forall\ m.
\end{equation*} 
\begin{proof}
The arguments to achieve the proof are very similar to the arguments used in the proof of Theorem \ref{th:saturato}. First of all let us consider the monomials 
\[
 F_\alpha = X^\alpha + \sum C_{\alpha\beta} X^\beta
\]
corresponding to the monomial basis $B_{\id{j}}$ of $\id{j}$ and the ideal $\id{h}(\id{j}) \subset k[C]$ of the stratum $\St_h(\id{j})$. 

In order to compute $\St_h(\id{j}_{\geqslant m})$, we have to consider again polynomials $F_\alpha$ as before if $\vert\alpha\vert \geqslant m,\ X^\alpha \in B_{\id{j}}$ and new polynomials $G_{\alpha\varepsilon}$ such that $\Lt(G_{\alpha\varepsilon}) = X^{\alpha+\varepsilon}$, $\forall\ X^\alpha \in B_{\id{j}},\ \vert\alpha\vert < m,$ and $\forall\ X^\varepsilon$ of degree $m-\vert\alpha\vert$, especially $X^\alpha X_n^{m-\vert\alpha\vert}$.
Then by the definition itself of $\mathtt{DegRevLex}$, the tail of $X^\alpha X_n^{m-\vert\alpha\vert}$ contains exactly the monomials in the tail of $X^\alpha$ multiplied by $X_n^{m-\vert\alpha\vert}$. So we can write
\[
G_{\alpha\varepsilon} = \begin{cases}
                         X^{\alpha+\varepsilon} + \sum E_{\alpha\delta}^\varepsilon X^\delta,& \forall\ X^\varepsilon \neq X_n^{m-\vert\alpha\vert},\\
                         X^\alpha X_n^{m-\vert\alpha\vert} + \sum C_{\alpha\beta} X^\beta X_n^{m-\vert\alpha\vert} = X_n^{m-\vert\alpha\vert} F_\alpha,&
                         \text{if } X^\varepsilon = X_n^{m-\vert\alpha\vert}\\
                        \end{cases}
\]
hence $\id{h}(\id{j}_{\geqslant m}) \subset k[C,E]$ (note that in the present case variables $D$ do not appear by construction).

By Theorem \ref{th:saturato} \emph{\ref{it:saturato_i})}, we know that all the variables $E$ can be eliminated. By the same reasoning used in the proof of Theorem \ref{th:saturato} \emph{\ref{it:saturato_ii})}, the ideal $\overline{\id{h}} = \id{h}(\id{j}_{\geqslant m}) \cap k[C]$ contains the $X$-coefficients of a set of S-polynomials corresponding to a set of the $S$-polynomials of the monomial basis of $\id{j}$: so $\St_h(\id{j}) \simeq \St_h(\id{j}_{\geqslant m})$.
\end{proof}
\end{prop}

\section{Gr\"obner strata and regularity} \label{sec:regularity}

In the present and following sections $k[X]$, $\prec$ and $\id j=(X^{\gamma_1}, \dots, X^{\gamma_t})$  will be as in the previous, but from now on we will consider only homogeneous ideals (with respect to the usual grading) and $T_i$ will be the complete homogeneous tail of $X^{\gamma_i}$ so that the only involved strata will be the homogeneous strata $\St_h(\id{j})$ introduced in Definition \ref{def:codaridotta} \emph{\ref{it:codaridotta_ii})}. Since every tail is fixed by $\prec$, we will simply denote ideals defining \Gr strata by $\id{h}(\id{j})$.

Let $p(z)$ be any admissible Hilbert polynomial for subschemes in $\PP^n$. Our goal is to show that the Hilbert scheme $\hilb_{p(z)}^n$ can be covered    by homogeneous strata of the type $\St_h(\id{j})$. In order to prove that, it is convenient to think of $\hilb_{p(z)}^n$ and $\St_h(\id{j})$ as schemes parameterizing the same kind of objects, namely homogeneous ideals in $k[X]$; as many ideals define the same subscheme  $Z \subset \PP^n$, the problem is to select a unique ideal in $k[X]$ for every subscheme $Z$. The most common choice is to associate to $Z$ the only homogeneous saturated ideal $I(Z)$ such that $Z=\Proj \left( k[X]/I(Z)\right)$; this point of view is that assumed for instance in \cite{NS} and in \cite{RT}, where homogeneous strata of saturated ideals are considered.

Here we prefer a different approach, that directly calls back to the explicit construction of Hilbert schemes
(see for instance \cite{B,HS,sern}).  

\begin{defin}\label{def:prt} Given an admissible Hilbert polynomial $p(z)$ for subschemes in $\PP^n$, we will denote by  $r$ the Gotzmann number of $p(z)$, that is the worst regularity of saturated ideals defining subschemes in $\hilbp$. Moreover we set: $M:=\binom{n+r}{n}$, $t:=M - p(r)$, $M_1:=\binom{n+r+1}{n} $ and $t_1:=M_1 - p(r+1)$. 
\end{defin}

Macaulay's Theorem states that $r$ is the regularity of the lexsegment ideal with Hilbert polynomial $p(z)$ (for the definition and the main properties of regularity and for  some consequences, we refer to \cite{Green-gin}).

As 
$Z=\Proj \left( k[X]/I(Z)\right)  = \Proj \left(  k[X]/I(Z)_{\geqslant r}\right)$, $Z$ can be uniquely  identified by the ideal $I(Z)_{\geqslant r}$, which is generated by $t$ linearly independent degree $r$ homogeneous polynomials $F_1, \dots, F_t$ or, more precisely, by the $t$-dimensional $k$-vector space $I(Z)_r$:  $\hilbp$ can be realized as a closed subscheme in the grassmannian of the $t$-dimensional vector spaces in $k[X]_r$.
A $t$-dimensional vector space in $k[X]_r$ gives a point in $\hilbp$ if and only if it generates an ideal $\id{i}$ having $p(z)$ as Hilbert polynomial.

\begin{notation}\label{not:hilb} From now on, $\id{i}\in \hilbp$ will mean that $\id{i}=I(Z)_{\geqslant r}$ for some closed subscheme $Z$ in $\PP^n$  with Hilbert polynomial $p(z)$. Equivalently we can say that $\id{i}\in \hilbp$ if and only if $\id{i}$ is an homogeneous ideal in $k[X]$ with Hilbert polynomial $p(z)$ (for the meaning of \lq\lq Hilbert polynomial of $\id{i} $\rq\rq\ see Notation \ref{not:hilb_scheme}) which  is generated in degree $r$, where $r$ is the Gotzmann number of $p(z)$.
\end{notation}

\begin{remark}\label{rem:saturation}
If $\id{i}  \in \hilbp$, then $\id{i}$ is $r$-regular and it has a free resolution of the type:
\begin{equation}\label{eq:resI_p-lineare}
0 \ \rightarrow\ k[X](-r-\lambda)^{n_\lambda}\ \rightarrow\ \ldots\ \rightarrow\ k[X](-r-1)^{n_1} \ \rightarrow\ k[X](-r)^{n_0} \ \rightarrow\ \id{i} \ \rightarrow\ 0 
\end{equation}
(\cite[Theorem 1.2]{EisGo}). Then we can find a set of generators for the first syzygies $\Syz (\id{i})$ in degree $r+1$.
\end{remark} 

If we take into consideration the  homogeneous Gr\"obner strata $\St_h(\id{j})$ and select  the    monomial ideal  $\id{j}$ in $\hilbp$, we  obtain the intended direct relation between Gr\"obner strata and Hilbert schemes. 

\begin{lemma}\label{th:sett} If $\id{j} \in \hilbp$, then (at least set-theoretically)  $\St_h(\id{j}) \subseteq \hilbp$. 
\begin{proof} Let $\id{i}$ be any ideal in $\St_h(\id{j})$. By hypothesis $\Lt(\id{i})=\id{j}$ and then $\id{i}$ and $\id{j}$ share the same Hilbert function. Therefore $\id{i}$ is generated in degree $r$ and has Hilbert polynomial $p(z)$ and then $\id{i} \in \hilbp$.
\end{proof}
\end{lemma}

Now we will see that the  set-theoretic inclusions are in fact algebraic maps and that for some ideals they are open injections. The crucial point is that the stratum structure (and so its injection in the Hilbert scheme) depends on the ideal $\id{j}$ and not on the the corresponding subscheme $Z=\Proj(k[X]/\id{j})$. This is not so surprising because the choice of the ideal fixes all the allowed deformations, but we want to stress this issue because in \cite{NS} the authors underestimated this fact and they made a wrong choice (proof of Corollary 4.4). 
In fact the stratum of the saturated lexicographic ideal $\id{L}$  with Hilbert polynomial $p(z)$ is not in general isomorphic to an open subset of $\hilbp$ (see \cite{RT} and Example \ref{rk:stratoLex}), whereas, as we will see,  the stratum of its truncation $\id{L}'=\id{L}_{\geqslant r}$  is an open subset of the Reeves-Stillman component of $\hilbp$.

\medskip

Let $\id{j}$ be a monomial ideal in $\hilbp$. As seen in \S\ \ref{sec:stratum_ideal} every ideal $\id{i}$ such that $\Lt(\id{i})=\id{j}$ has a (unique) reduced Gr\"obner basis $\{f_1, \dots, f_t\}$ where $f_i$ is as in Definition \ref{def:codaridotta} \emph{\ref{it:codaridotta_ii})}. Not every ideal generated by $t$ polynomials of such a type has $\id{j}$ as initial ideal. In order to obtain equations for $\St_h(\id{j})$ we  consider the coefficients $c_{i\alpha}$ appearing in the $f_i$ as new variables; more precisely let  $C=\{ C_{i\alpha},\ i=1, \dots, t, \ X^{\alpha}\in k[X]_r \setminus \id{j}_r \text{ and } X^\alpha \prec X^{\gamma_i} \}$ be new variables and consider $t$ polynomials in $k[X,C]$ of the following type:
\begin{equation}\label{eq:costF}
  F_i = X^{\gamma_i} + \sum_{X^{\alpha}\in T_i}  {C}_{i\alpha} X^\alpha  
\end{equation}
where $T_i= T_{\gamma_i}\cap k[X]_r$ (Definition \ref{def:tail}). We obtain the ideal $\id{h}(\id{j})$ of $\St_h(\id{j})$ collecting the $X$-coefficients of some complete reduction with respect to $F_1, \dots, F_t$ of all the $S$-polynomials $S(F_i,F_j)$, corresponding to a set of generators for $\Syz(\id{j})$ (see Theorem \ref{th:cor-unica} and Proposition \ref{th:unica} \emph{\ref{it:unica_4})}).

\begin{prop}\label{th:matriceA} In the above notation, let $\id{j}$ be a monomial ideal in $\hilbp$ and let $A$ be the $t(n+1)\times M_1$ matrix whose entries are the $X$-coefficients of $X_jF_i$, for all $j = 0,\dots,n$ and $i = 1,\dots,t$.

Then the ideal $\id{h}(\id{j})$ of the homogeneous stratum $\St_h(\id{j})$ is generated by the $(t_1+1) \times (t_1 + 1)$ minors of $A$.
\begin{proof} By abuse of notation we write in the same way a polynomial and the rows of its $X$-coefficients.   As in Definition \ref{def:procedura} we consider a term order on $\T_{X,C}$ which is an elimination order of the variables $X$ and coincides with the fixed term ordering $\prec$ on $\T_X$. It is quite evident by elementary arguments of linear algebra, that the ideal $\id{a}\subseteq k[C]$, generated by all $(t_1+1) \times (t_1 + 1)$ minors, does not change if we perform some row reduction on $A$. Let $\mathcal{P}$ be a set of $t_1$ rows whose leading terms are a basis of $\id{j}_{r+1}$. If $X_hF_i \notin \mathcal P$, then it has the same leading term than one in $\mathcal{P}$, say $X_kF_j$; we can substitute $X_hF_i$ with $X_hF_i-X_kF_j$. In this way the rows not in $\mathcal{P}$ become precisely all the $S$-polynomials $S(F_i,F_j)$ that have $X$-degree $r+1$.

\medskip

At the end of this sequence of row reductions, we can write the matrix as follows:
\begin{equation}\label{eq:matriceA1lineare}
\left(\begin{array}{c|c}
D & E \\
\hline
S & L
\end{array}\right)
\end{equation}
where $D$ is a $t_1 \times t_1$ upper-triangular matrix with 1's along the main diagonal, whose rows correspond to $\mathcal{P}$ and whose columns correspond to monomials in $\id{j}_{r+1}$.

\medskip

Using rows in $\mathcal{P}$, we now perform a sequence of rows reductions on the following ones, in order to annihilate all the coefficients of monomials in $\id{j}_{r+1}$, that is the entries of the submatrix $S$: if $a(C)$ is the first non-zero entry in a row not in $\mathcal{P}$ and its column corresponds to the monomial $X^{\gamma} \in \id{j}_{r+1}$, we add to this row $-a(C)X_kF_j$, where $X_kF_j \in \mathcal{P}$ and $\Lt(X_kF_j)=X^{\gamma}$. This is nothing else than a step of reduction with respect to $\{F_1, \dots, F_t\}$. At the end of this second turn of rows reductions, we can write the matrix as follows:
\begin{equation}\label{eq:matriceA1strato}
\left( \begin{array}{c|c}
D & E \\
\hline
0 & R 
\end{array}\right) 
\end{equation} 
where the rows in $(D\ |\ E)$ are unchanged whereas the rows in $(0\ |\ R)$ are the $X$-coefficients of     complete reductions of $S$-polynomials in $X$-degree $r+1$. Then $\id{a}$ is generated by the entries of $R$ and so  $\id{a} \subset \id{h}(\id{j})$.

We can see that this inclusion is in fact an equality taking in mind Remark \ref{rem:saturation} and Proposition \ref{th:unica} \emph{\ref{it:unica_4})}: the first one says that $\Syz(\id{j})$ is generated in degree $r+1$ and the second one that in this case $\id{h}(\id{j})$ is generated by the $X$-coefficients of complete reductions of the $S$-polynomials $S(F_i,F_j)$ of $X$-degree $r+1$.
\end{proof}
\end{prop}

The following corollary just express in an explicit way two properties contained in the proof of Proposition \ref{th:matriceA}.

\begin{cor}\label{cor:matriceA} In the above notation:
\begin{itemize}
\item the ideal $\id{h}(\id{j})$ is generated by the entries of the submatrix $R$ in \eqref{eq:matriceA1strato};
\item the vector space $L(\id{j})$ is generated by the entries of the submatrix $L$ in \eqref{eq:matriceA1lineare}.
\end{itemize}
\end{cor}

As already said in the Remark \ref{rk:termordersignificance}, this theorem shows one more time that Gr\"obner strata equations are substantially independent of the term ordering, that sets only which monomials can appear in the tails $T_i$.

\section{\Gr strata that are open subsets of an Hilbert scheme}
\label{sec:opensubset}  

In the present section we will prove that every homogeneous Gr\"obner stratum $\St_h(\id{j})$, where $\id{j} \in \hilbp$, can be naturally identified with a locally closed subscheme of $\hilbp$ and that it is an open subset of $\hilbp$ if $\id{j}$ is generated by the first $t$ monomials in $k[X]$ with respect to the fixed term ordering $\prec$. As a consequence we obtain the main results of the paper about the rationality of some components of $\hilbp$.

For the meaning of $p(z)$, $r$, $t$, $M$, $t_1$, $M_1$ and $\id{i} \in \hilbp$ we refer to Definition \ref{def:prt} and Notation \ref{not:hilb}.

First of all we recall how equations defining $\hilbp$ are usually obtained (see for instance \cite{B,HS}). Every ideal $\id{i} \in \hilbp$ is generated by the $t$-dimensional vector space $\id{i}_r$. On the other hand, thanks to Gotzmann's Persistence Theorem  (see for instance \cite[Theorem 3.8]{Green-gin}),  a $t$-dimensional vector space $V \subset  k[X]_r$ generates an ideal $\id{i} \in \hilbp$ if and only if $\dim_k \langle X_0V, \dots, X_nV\rangle =t_1$.

Therefore $\hilbp$ can be thought as the subscheme of the grassmannian $\mathbb{G}(t,M)$ defined by the previous condition. Moreover by the Pl\"{u}cker embedding of the grassmannian in a projective space $\PP^q$, $\hilbp$ becomes a closed subscheme (not necessarily irreducible and reduced) of $\PP^q$.

Here we are not interested in finding explicit equations for $\hilbp$ in $\PP^q$, but only equations defining each open subset $U \cap \hilbp$, where $U$ is the open subset of $\mathbb{G}(t,M)$ given by a non-vanishing Pl\"{u}cker coordinate.  

\begin{defin}\label{plu} Thinking of $k[X]_r$ as the vector space generated by its monomials, we can identify every Pl\"{u}cker coordinate with a suitable monomial ideal $\id{j}$ generated by $t$ monomials of degree $r$. We will denote by $U_{\id{j}}$ and $\HH_{\id{j}}$ respectively the open subsets of $\mathbb{G}(t,M)$ and of $\hilbp$ where the Pl\"{u}cker coordinate corresponding to $\id{j}$ does not vanish. 
\end{defin}

In a natural way $U_{\id{j}}$ is isomorphic to the affine space  $\mathbb{A}^{t(M-t)}$. In fact,  if $\id{j}=(X^{\gamma_1}, \dots, X^{\gamma_t})$, every point in $U_{\id{j}}$ is uniquely identified by   the reduced, ordered set of generators $\langle g _1, \dots,g _t \rangle $  of the type  $g_i = X^{\gamma_i}+ \sum c_{i\alpha} X^\alpha$, where  $c_{i\alpha}\in k$ and  $X^{\alpha}$ is any monomial in $ k[X]_r\setminus \idj$. Then   we consider on $\mathbb{A}^{t(M-t)}$ the coordinates   $C_{i\alpha}$. Note that each $C_{i\alpha}$  naturally corresponds to the  Pl\"{u}cker coordinate $\idj'=( X^{\gamma_1}, \dots,X^{\gamma_{i-1}},X^{\alpha} ,X^{\gamma_{i+1}},\dots  X^{\gamma_t} )$ (but of course not all the \Pl coordinates are of this type).

\medskip

Now we can mimic the construction of Gr\"obner strata and obtain the  defining  ideal  of $\HH_\idj$ as a subscheme of  $\mathbb{A}^{t(M-t)}$. Let us consider the set of variables $\overline{C}=\{ C_{i\alpha}, i=1, \dots, t, \ X^{\alpha}\in k[X]_r \setminus \id{j}\}$ and $t$ polynomials $G_1, \dots, G_t$ in $k[X,\overline{C}]$ of the type:
\begin{equation}\label{eq:costG} 
G_i = X^{\gamma_i}+\sum C_{i\alpha} X^\alpha   
\end{equation}
and let $B$ be the $(n+1)t\times M_1$ matrix whose entries are the $X$-coefficients of the polynomials $X_jG_i$. Then consider the ideal  $\id{b}(\id{j}) \subset k[\overline{C}]$  generated by the $(t_1+1)\times (t_1+1)$ minors in $B$. 

\begin{prop}\label{th:rkB}  $\id{b}(\idj)$ is the ideal of $\HH_\idj$ as a closed subscheme of $\mathbb{A}^{t(M-t)}$.
\end{prop}
\begin{proof}
Every ideal $\id{i} \in U_\idj$ can be obtained from $(G_1, \dots, G_t)$ specializing (in a unique way) the variables $C_{i\alpha}$ to $c_{i\alpha}\in k$. Obviously not all the specializations give ideals $\id{i} \in \HH_\idj$, that is  with Hilbert polynomial $p(z)$ (more precisely, such that $k[X_0, \dots, X_n]/\id{i}$  has Hilbert polynomial $p(z)$), because we have to ask both $\dim_k(\id{i}_r)=t$ and $\dim_k \id{i}_{r+1} = t_1 $: thanks to Gotzmann's persistence we know that these two necessary conditions are also sufficient. 

In the open subset $U_{\id{j}}$ the first condition always holds and the rank of every specialization of $B$ is $\geqslant t_1$ by Macaulay estimate of the growth of ideals (see \cite[Section 3]{Green-gin} or \cite[Corollary 5.5.28]{KrRob2}). Therefore $\HH_\idj$ is given by the condition $\rk(B)\leqslant t_1$.
\end{proof}

We  can order   the set of \Pl coordinates in the following way. We write the  $t$ monomials  corresponding to each \Pl coordinate  in decreasing order with respect to $\prec$; if $\idj_1=(\ X^{\alpha_1}\succ  \dots \succ X^{\alpha_t})$ and $\idj_2=(\ X^{\beta_1}\succ  \dots \succ X^{\beta_t})$,  then   $\idj_1\succ \idj_2$  if $\ X^{\alpha_i}= X^{\beta_i}$ for every $i$ lower than some $s$ and $\ X^{\alpha_s}\succ X^{\beta_s}$.

\medskip

It is now easy to compare, for the same monomial ideal $\id{j} \in \hilbp$, the Gr\"obner stratum $\St_h(\id{j})$ and the open subset $\HH_{\id{j}}$.  We underline that for our purpose it will be sufficient to consider the open subsets $\HH_{\id{j}}$ corresponding to monomial ideals  $\id{j} \in \hilbp$,  because (scheme-theoretically) they  cover  $\hilbp$. In fact, if $\id{i} \in \hilbp$, then also $\Lt(\id{i}) \in \hilbp$ and so  $\id{i} \in \HH_{\textnormal{\scriptsize LT}(\id{i})}$.

\begin{theorem}\label{th:aperto} Let $p(z)$ be any admissible Hilbert polynomial in $\PP^n$ with Gotzmann number $r$.    Let us fix  any term ordering $\prec$ on $\T_X$. 
\begin{enumerate}[i)]
\item\label{it:aperto_i}  If $\id{j}$ is a monomial ideal in $ \hilbp$, then $\St_h(\id{j})$ is naturally isomorphic to the locally closed subscheme of $\hilbp$ given by the conditions that the Pl\"{u}cker coordinate corresponding to $\id{j}$ does not vanish and the preceding  ones vanish.
\item\label{it:aperto_ii} For every  isolated, irreducible component  $H$ of $\hilbp$,  there is a monomial ideal $\id{j} \in \hilbp$ such that an irreducible  component of $\Supp \St_h(\id{j})$ is an open subset of $\Supp H$. Then 
$\Supp H$ has an open subset   which is a homogeneous affine variety with respect to a non-standard positive grading.
\item\label{it:aperto_iii} Every smooth irreducible component $H$ of $\hilbp$ is rational. 
The same holds for every smooth, irreducible component of $\Supp \hilbp$.
\end{enumerate}
\begin{proof} 
\emph{\ref{it:aperto_i})} We obtain the two affine varieties $\St_h(\id{j})$ and $\HH_{\id{j}}$ in a quite similar way (for $\St_h(\id{j})$ see Proposition \ref{th:matriceA} and for $\HH_\idj$ see Proposition \ref{th:rkB}). The only difference comes from the definition of the set of polynomials $F_1, \dots, F_t$ given in \eqref{eq:costF}, leading to equations for $\St_h(\id{j})$, and the set of polynomials $G_1, \dots, G_t$ given in \eqref{eq:costG}, leading to equations for $\HH_{\id{j}}$: in $G_i$ the sum is over all the degree $r$ monomials $X^{\alpha}\notin \id{j}$ whereas in $F_i$ we also assume that $X^{\alpha} \prec \Lt(F_i)$. Therefore we can think of $\St_h(\id{j})$ as the affine subscheme defined by the ideal $\overline{\id{h}(\idj)}$ in the ring $k[X,\overline{C}]$, where $\overline{C}=\{C_{i\alpha} \ | \ i=0, \dots, n, \  \ X^{\alpha} \in k[X]_r \setminus \id{j} \}$ generated by $\id{h}(\id{j})$ and by $\big(C_{i\alpha}\  |\  X^{\alpha} \succ \Lt(F_i)\big)$, namely   $\overline{\id{h}(\idj)}=\id{h}(\idj)k[\overline{C}]+( \overline{C}\setminus C)$. Now we can  conclude because  the Pl\"{u}cker coordinates higher than $\id{j}$ vanish if and only if all the $C_{i\alpha}$ such that $ X^{\alpha} \succ \Lt(F_i)$ vanish.

\emph{\ref{it:aperto_ii})} As $\id{j}$ varies among the finite set of the monomial ideals in $\hilbp$, the \Gr strata $\St_h(\id{j})$ give a set theoretical covering of $\hilbp$ by locally closed subschemes. Then there is a suitable ideal $\id{j}$ such that an irreducible component of $\Supp \St_h(\id{j})$ is an open subset of $H$.
We have seen in the previous sections that  $\St_h(\id{j})$ has a structure of homogeneous affine scheme with respect to a non-standard positive grading $\lambda$. Then also its support and the irreducible components of the support are homogeneous (see \cite[Section IV.3.3]{Bourbaki} and \cite[Corollary 2.7]{FR}).

\emph{\ref{it:aperto_iii})} If $H$ is a smooth, irreducible component of either $\hilbp$ or $\Supp \hilbp$, then it is also reduced. Thanks to the previous item we know that an open subset of $H$ is an affine homogeneous variety with respect to a positive grading. Moreover this open subset is also smooth and so  it is isomorphic to an affine space, by 
 Corollary \ref{th:liscio}. 
\end{proof}
\end{theorem}

\begin{remark} Let $\id{j}$ be a monomial ideal in $\hilbp$ and let $\id{b}(\id{j}) \subset k[\overline{C}]$ the ideal of $\HH_{\id{j}}$. It is possible to define a grading $\lambda'$ on $k[\overline{C}]$ such that $\id{b}(\id{j})$ becomes homogeneous, by the analogous definition: $\lambda'(C_{i\alpha})=\gamma_i-\alpha$ if $C_{i\alpha}$ appears in $G_i$ \eqref{eq:costG}. However  this grading $\lambda'$ is not necessarily   positive and so it gives less interesting consequences.

If an irreducible  component $H$ of $\hilbp$ is also reduced, Theorem \ref{th:aperto} insures that there is an open subset of $H$ which has the structure of homogeneous variety with respect to a positive grading induced  from that of a suitable \Gr stratum $\St_h(\idj)$.

On the other hand, in the case of  a non-reduced component we only know that the support of a suitable   open subset  is   homogeneous with respect to a positive grading, but this does not imply that the open subset itself is homogeneous. 
\end{remark}

Now we consider a special case in which we obtain a positive grading on an open subset of an irreducible component of $\hilbp$, even if not reduced.

\begin{defin}\label{def:segment} 
Given any term order $\prec$ in $\T_X$, a $(m,\prec)$-\emph{segment} is a subset $S$ of $k[X]_m$ containing  the first $\vert S \vert$ monomials of degree $m$ with respect to $\prec$, namely such that:
\[
\forall\ X^\beta \in k[X]_m,\ \forall\ X^\gamma \in S:\  X^{\beta } \succ X^\gamma \Rightarrow X^\beta \in S.
\]
An $(m,\prec)$\emph{-segment ideal }is a monomial ideal $\id{j}$ which is generated by a $(m,\prec)$-segment.
\end{defin}

If $\id{L}$ is the saturated lexsegment ideal, then for every $m \geqslant r$ (that is for every $m$ higher than the regularity of $\id{L}$), the ideal $\id{L}_{\geqslant m}$ is a $(m,\mathtt{Lex})$-segment ideal. This property does not hold in general if the term ordering is not $\mathtt{Lex}$, so that $\id{j}_{\geqslant m}$ could be a  $(m,\prec)$-segment ideal and $\id{j}_{\geqslant m+1}$ could not be a $(m+1,\prec)$-segment ideal. A trivial case is for instance that of the ideal $(X_0)$ in $k[X_0,X_1,X_2]$ which is $(1,\mathtt{DegRevLex})$-segment ideal, whereas $(X_0)_{\geqslant 2}=(X_0^2,X_0 X_1,X_0 X_2)$ is not a $(2,\mathtt{DegRevLex})$-segment ideal, because it contains $X_0 X_2$ and does not contain $X_1^2$. 

The  definition of $(m,\prec)$-segment ideal is not equivalent, but it is very  close to that of \emph{extremal ideal} given in \cite{sherman-2007}.

\begin{cor}\label{cor:vuoto} Let $\id{j}$ be $(r,\prec)$-segment ideal in the grassmannian $\mathbb{G}(t,M)$. 

If $\id{j}$ does not belong to $\hilbp$, then the open subset $\HH_{\id{j}}$ is empty.
\begin{proof} Any point $\id{i} \in \HH_{\id{j}}$ should belong to the Gr\"obner stratum $\St_{h}(\id{j})$, that is it should share the same Hilbert polynomial of $\id{j}$, which is not $p(z)$.
\end{proof}
\end{cor}

The first of the following examples  highlights both that  Theorem \ref{th:aperto} does not hold for a monomial ideal $\idj $ that belongs to $\mathbb{G}(t,M)$   but  not to $\hilbp$ and that Corollary \ref{cor:vuoto} does not hold  for a monomial ideal $\idj$ in  $\mathbb{G}(t,M)$ which is not a segment. Moreover Example \ref{ex:got} presents a concrete case of empty $\HH_{\id{j}}$ as discussed in the previous corollary.

\begin{example} Let us consider the constant Hilbert polynomial $p(z)=2$ on $\PP^2$. As well known, $\hilb^2_2$ is irreducible of dimension 4. The monomial ideal $\idj =(X_0^2,X_0X_1,X_1^2,X_2^2)$ is generated by $4$ monomials of degree $2$, but does not belong to $\hilb^2_2$ because its radical is the irrelevant maximal ideal. However, $\mathcal{H}_\idj$ is non-empty because it contains for instance all the reduced subschemes given by  couples  of points $P[1:a:b],\ Q[1:a':b']\in \PP^2$ such that $ab'\neq a'b$. By the way, $\St_h(\id{j})$ cannot have any common point with $\hilb^2_2$.
\end{example}

\begin{example}\label{ex:got}  In the example presented at the end of  \S\ \ref{sec:example} the complete list of Borel ideals in $k[X_0,X_1,X_2,X_3]$ with Hilbert polynomial $p(z)=4z$ is presented. None of them is the $(6,\texttt{DegRevLex}$)-segment ideal $\id{j}$ containing all the $t=60$ degree $6$ monomials except the $p(6)=24$ lowest with respect to the term ordering $\texttt{DegRevLex}$. As $r=6$ is the Gotzmann number of $4z$ and $M=84$, then   $\id{j}$ belongs  to $\mathbb{G}(t,M)$ and does not to $\hilb_{4z}^3$ ($\id{j}$ has constant Hilbert polynomial equal to 24). Hence $\HH_{\id{j}}$ is empty.
\end{example}

\begin{cor}\label{th:segment01} Let $p(z)$ any admissible Hilbert polynomial in $\PP^n$ with Gotzmann number $r$ and let $H$ be an isolated, irreducible component of $\hilbp$.  

If $H$ contains a point corresponding to an $(r,\prec)$-segment ideal $\id{j}\in \hilbp$ with respect to some term ordering $\prec$ on $\T_X$, then $\St_h(\id{j})$ is an open subset of $H$, so that 
$H$ has an open subset   which is an homogeneous affine variety with respect to a non-standard positive grading.

\begin{proof} If $\id{j}$ is a $(r,\prec)$-segment ideal, then there are no Pl\"{u}cker coordinates preceding    that corresponding to $\id{j}$. Thus $\St_h(\idj)\cong \HH_\idj$ (see Theorem \ref{th:aperto}) and so $\HH_\idj$ is an affine  homogeneous scheme with respect to a positive grading.   
\end{proof}
\end{cor}

\begin{cor}\label{th:segment}
Let  $\id{j} \in \hilbp$ be $(r,\prec)$-segment ideal and  let $H$ be an irreducible component of $\hilbp$ containing $\id{j}$. If either of the following condition holds: 
\begin{enumerate}[i)]
\item  $\St_h(\id{j})$ is an affine space,
\item  $\id{j}$ is a smooth point of   $\St_h(\id{j})$ ,
\item  $\id{j}$ is a smooth point of  $\hilbp$,
\end{enumerate}
then $H$ is rational.
\end{cor}

\begin{proof}
Straightforward consequence of the previous result and of Corollary \ref{th:liscio}.
\end{proof}

\section{Gr\"obner stratum of the lexsegment ideal} \label{sec:lexsegment} 

In this paragraph the term ordering $\prec$ will be the lexicographic term ordering \texttt{Lex}.

As a first application to the results obtained in \S\ \ref{sec:opensubset}, we take into consideration the  \emph{lexicographic ideal} $\id{L}$. For every admissible Hilbert function $p(z)$ on $\PP^n$, $\hilbp$ contains the ideal generated by the first $t$ monomials in degree $r$ with respect to the term ordering \texttt{Lex}. In the paper \cite{RS} it is proved that the point of  $\hilbp$ (usually called \emph{lexicographic point}) corresponding to the subscheme $\Proj k[X]/\id{L}$ is smooth and Reeves and Stillman get the proof by a computation of the Zariski tangent space dimension. The only component of $\hilbp$ containing the lexicographic point is usually denoted by $H_{RS}$.
 As a consequence of the quoted result by Reeves and Stillman and  of Corollary \ref{th:segment}, we then obtain:

\begin{cor}\label{th:lexlex}
The Reeves and Stillman component $H_{RS}$ of $\hilbp$ is rational.
\end{cor}

 However we prefer to present here a new, self-contained proof, in order to explain how our technique can be used as a theoretical, as well as a computational, tool.

First of all, we recall briefly the notation used in \cite{RS}. Moving from \cite{Mac}, Reeves and Stillman work with lexicographic saturated ideals of the type:
\[
L(a_0,\dots,a_{n-1}) = (X_0^{a_{n-1}+1},X_0^{a_{n-1}}X_1^{a_{n-2}+1},\dots,X_0^{a_{n-1}}\cdots X_{n-2}^{a_1}X_{n-1}^{a_0}).
\]
Since we are going to prove the same result, we will assume a quite similar notation, but not the same because  in this paper we consider ideals generated in degree $r$ instead of saturated ideals.

\begin{notation}
We refer with $\id{L}(a_0,\ldots,a_{n-1})$ to the lexsegment ideal generated by all monomials of degree $r = \sum a_j$ that precede (greater than or equal to) the monomial $X_0^{a_{n-1}}\cdots X_{n-1}^{a_0}$ in the \texttt{Lex} term ordering:
\[
 \id{L}(a_0,\ldots,a_{n-1}) = (X_0^r,X_0^{r-1}X_{1},\ldots,X_0^{a_{n-1}}\cdots X_{n-1}^{a_0}).
\]
Note that $r$ is precisely the Gotzmann number of the Hilbert polynomial of $\id{L}$ (more precisely of $k[X]/\id{L}$).
\end{notation}

\begin{theorem}\label{th:lexstratumopenrational}
The homogeneous \Gr stratum $\St_h(\id{L}(a_0,\ldots,a_{n-1}))$ of the lexicographic ideal  $\id{L}(a_0,\ldots,a_{n-1})\in \hilbp$ is  isomorphic to an affine space. Therefore the  component $H_{RS}$ of $\hilbp$ is rational.
\end{theorem}

\begin{proof}
Thanks to Corollary \ref{th:segment} we obtain the complete statement proving that the homogeneous \Gr stratum $\St_h(\id{L}(a_0,\dots,a_{n-1}))$ is an affine space, that is showing that a same number is both a lower-bound for its dimension and an upper-bound for its embedding dimension; the first part corresponds to Theorem 4.1 of \cite{RS} (here in terms of initial ideals) and the second one corresponds to Theorem 3.3 of \cite{RS}.

We proceed by induction on the number $n$ of variables and on the Gotzmann number  $r$.

In order to obtain an upper-bound for the embedding dimension we look for a maximal set of eliminable variables $C'\subset C$, using Criterion \ref{criterio}. If $X^{\alpha_1} \succ \dots \succ X^{\alpha_m}$ is the monomial basis of the saturation $\id{l}$ of $\id{L}(a_0,\dots, a_{n-1})$, then we can assume that the polynomials $F_1, \dots, F_t\in k[X,C]$ (that we use in order to construct $\St_h(\id{L}(a_0,\dots, a_{n-1}))$: see Definition \ref{def:procedura}) are ordered so that $\Lt(F_i)=X^{\alpha_i}X_n^{r-\vert \alpha_i \vert}$ for $i=1,\dots, m$. Thanks to Theorem \ref{th:saturato} we can start the construction of $C'$, putting inside all the variables appearing in $F_j$ for every $j>m$.

We divide the proof in 3 steps.

\medskip

\textbf{Step 1}  The zero-dimensional case: $ \St_h(\id{L}(a_0,0, \dots, 0))\simeq \Af^{na_0}$.

\medskip

\noindent \textbf{Claim 1i:}  $\dim \St_h(\id{L}(a_0,0,\dots,0))\geqslant na_0$.

\noindent  Let us denote $\id{L}(a_0,0, \dots, 0)$ by $\id{L}$. The zero-dimensional scheme $Z$ of $a_0$ general points in $\PP^n$ has Gotzmann number $a_0$ and Hilbert polynomial $p(z)=a_0$. Moreover $\Lt\big(I(Z)_{\geqslant a_0}\big)\supseteq \id{L}$, because for every monomial $X^{\gamma} \succeq X_1^{a_1}$ we can find some homogeneous polynomial of the type $X^{\gamma}-\sum_{j=1}^{a_0}b_j X_{n-1}^{a_0-j} X_n^j$ vanishing in the $a_0$ points of $Z$: we can find the $b_j$'s solving a $a_0\times a_0$ linear system with a Vandermonde associated matrix. As both $\Lt\big(I(Z)_{\geqslant a_0}\big)$ and $\id{L}$ are generated in degree $a_0$, they coincide; so $I(Z)_{\geqslant r} \in \St_h(\id{L})$ and we conclude since we can choose $Z$ in a family of dimension $na_0$.

\medskip

\noindent \textbf{Claim  1ii:} $\ed\, \St_h(\id{L}(a_0,0, \dots, 0))\leqslant na_0$. 

\medskip

\noindent The saturation of $\id{L}$ is the ideal $(X_0,X_{1},\ldots,X_{n-2},X_{n-1}^{a_0})$, which is generated by $n$ monomials; moreover there are only $a_0$ monomials of degree $a_0$ not contained in $\id{L})$: Corollary \ref{th:saturato} leads to the conclusion.

\medskip

\textbf{Step 2:}  If $\St_h(\id{L}(0,a_1,\dots,a_{n-1})) \simeq \Af^K$ then $ \St_h(\id{L}(a_0,a_1, \dots, a_{n-1}))\simeq \Af^{K+na_0}$.

\medskip

\noindent \textbf{Claim 2i:} $\dim \St_h(\id{L}(a_0,a_1, \dots, a_{n-1})) \geqslant \dim  \St_h(\id{L}(0,a_1,\dots,a_{n-1}))+na_0=K+na_0$.

\medskip

\noindent Let us denote $\id{L}(a_0,a_1, \dots, a_{n-1})$ by $\id{L}$ and $\id{L}(0,a_1, \dots, a_{n-1})$  by $\id{L}_1$.
Let $Y$ be any closed subscheme in $\PP^n$ such that $I(Y)_{\geqslant r}\in \St_h(\id{L}_1)$ and consider the set $Z$ of $a_0$ points in $\PP^n$. If we choose the $a_0$ points in $Z$ general enough, then $I(Z\cup Y)=I(Z)\cdot I(Y)$. Then we conclude thanks to the previous step, as $\Lt(I(Z))=\id{L}(a_0,0,\dots, 0)$ and $\id{L}=\id{L}_1\cdot \id{L}(a_0,0, \dots,0)$.

\medskip

\noindent \textbf{Claim 2ii:} $\ed\, \St_h(\id{L}(a_0,a_1, \dots, a_{n-1})) \leqslant \ed \, \St_h(\id{L}(0,a_1, \dots, a_{n-1}))+  na_0=K+na_0$.

\medskip

First of all, let us consider all the polynomials $F_i$ such that $X_n^{r-a_0}\mid  \Lt(F_i)$ and  the set of variables $C_{i\beta}$ appearing in them such that $X^{\beta}=X^{\beta_1}X_n^{r-a_1} $ for some monomial $X^{\beta_1}\notin \id{L}_1$: a multiple of  $X^\beta$ belongs to $\id{L}$ if and only the corresponding multiple of $X^{\beta_1}$  belongs to $\id{L}_1$. Then $F_i=X_n^{r-a_0}F_{1i}+\dots $, where the $F_{i1}$'s are the polynomials  that appear in the definition of $\St_h(\id{L}_1)$. Using the $S$-polynomials involving couples of such polynomials we see that $L(\id{L}_1)\subseteq L(\id{L})$; thus all the variables $C_{i\beta}$ of this type are eliminable, except at most $K=\ed\, \id{L}_1$ of them.

Moreover, for every $i \leqslant n$ there are $a_0$ variables $C_{i\beta}$ such that $X^{\beta}\notin  \id{L}$, $X^{\beta}\in  \id{L}_1$: the are  $X_0^{a_{n-1}}\cdots X_{n-2}^{a_1}X_{n-1}^{a_0-j}X_n^j$, $j=1, \dots,a_0$.

If we specialize to 0 all the variables of the two above considered types, the embedding dimension drops at most by $\ed\, \St_h(\id{L}_1)+na_0=K+na_0$. 

Now it will be sufficient to verify that all the remaining variables $C_{i\beta}$ are eliminable, using Criterion \ref{criterio}.

Assume that $ X^{\beta}\prec X_0^{a_{n-1}}\cdots X_{n-2}^{a_1}$ and $X_n^{r-a_0}\nmid X^{\beta}$.
\begin{itemize}
\item If $i>n$, all the variables are eliminable using those appearing in $F_1, \dots, F_n$, thanks to Corollary \ref{th:saturato}.
\item If $i<n$, using $S(F_i,F_j)$, where $\Lt(F_j)=X^{\alpha_i}X_{n-1}^{r-\vert \alpha_i \vert}$, we see that $C_{i\beta} \in L(\id{L})$.  
\item If $i=n$,  using $S(F_n,F_{n-1})=X_{n-2}X_n^{a_0-1}F_n-X_{n-1}^{a_0}F_{n-1}$, we see that $C_{n\beta} \in L(\id{L})$ (note that by the previous idem  $C_{n-1,\beta'} \in L(\id{L})$).	
\end{itemize}

\medskip 

\textbf{Step 3:} If $ \id{L}(0,a_1,\dots,a_{n-1}) \simeq \Af^{K_1}$  then $ \id{L}(0,a_1,\dots, a_d)\simeq \Af^{K_2}$ where $d$ is the maximal index $< n$ such that $a_d\neq 0$ and $ K_2=K_1+(n-d)(d+1)+\binom{a_{n-1}+n}{n}-1$ (or $K_2=K_1+\binom{a_{n-1}+n}{n}-1$ if $d$ does not exist).

\medskip

\noindent Here we compare the ideal $\id{L}= \id{L}(0,a_1,\dots,a_{n-1})$ in $k[X]$ and the ideal $\id{L}_1= \id{L}(0,a_1,\dots, a_d)$ in $k[X_0, \dots, X_d]$.
Observe that both $\id{l}:=\sat(\id{L})$ and $\id{l}_1:=\sat(\id{L}_1)$ fulfill the hypothesis of Theorem \ref{th:saturato} \emph{\ref{it:saturato_ii})} (see also Example \ref{rk:stratoLex}); then it holds $\St_h(\id{L})\simeq \St_h(\id{l})$ and $\St_h(\id{L}_1)\simeq \St_h(\id{l}_1)$. The statement for the saturated ideals $\id{l}$ and $\id{l}_1$ is proved using the same technique as above in \cite[Proposition 4.5]{RT}.
\end{proof}

\section{Algorithms and Examples} \label{sec:example}

In this final section we exhibit a pseudo-code description of the algorithms to compute the ideal of Gr\"obner stratum (Algorithm \ref{alg:idSth}) and its embedding dimension (Algorithm \ref{alg:ed}) and then we apply our technique to the Hilbert scheme $\hilb_{4z}^3$.

\subsection*{Algorithms}
The following two algorithms are mainly based on Proposition \ref{th:matriceA} and Corollary \ref{cor:matriceA}.
\begin{algorithm}[H]
\caption{\label{alg:idSth}Computing the ideal $\id{h}(\id{j})$ of $\St_h(\id{j})$}
\begin{algorithmic}[1] 
\Procedure{Gr\"obnerStratum}{$\id{j}$,$\prec$}
 \State Compute $B_\id{j}$ minimal monomial basis of $\id{j}$;
 \State $\mathcal{G} \leftarrow \emptyset$;
 \ForAll{$X^\alpha \in B_\id{j}$}
   \State Compute $F_\alpha = X^\alpha + \sum C_{\alpha\beta} X^\beta,\quad X^\beta \notin \id{j},\ X^\alpha \succ X^\beta,\ \vert\alpha\vert = \vert\beta\vert$;
   \State $\mathcal{G} \leftarrow \mathcal{G} \cup \{F_\alpha\}$;
 \EndFor
 \State Compute $\Syz(\id{j})$ basis of the syzygies of $B_\id{j}$;
 \State $\mathcal{C} \leftarrow \emptyset$;
 \ForAll{$(\ldots,X^{\delta_i},\ldots,X^{\delta_j},\ldots) \in \Syz(\id{j})$}
 \State $S \leftarrow X^{\delta_i}F_{\alpha_i} - X^{\delta_j}F_{\alpha_j}$;
 \State Compute the reduction $S_{\textnormal{red}}$ of $S$ w.r.t. the set of polynomials $\mathcal{G}$;
 \State Collect the set $\mathcal{R}$ of the $X$-coefficients of $S_{\textnormal{red}}$;
 \State $\mathcal{C} \leftarrow \mathcal{C} \cup \mathcal{R}$;
 \EndFor
 \State $\id{h} \leftarrow 0$;
 \While{$\mathcal{C} \neq \emptyset$}
   \State $g \leftarrow \min(\mathcal{C})$; \Comment{The minimun w.r.t. the $\lambda$-grading induced by $\prec$ on $k[C]$}
    \If{$L(g) = 0$} \Comment{$g$ has no linear part}
       \State $\id{h} \leftarrow \id{h} + g$;
    \Else
       \State Use $g$ to eliminate one variable from all the polynomials in $\mathcal{C}$;
    \EndIf
   \State $\mathcal{C} \leftarrow \mathcal{C} \setminus \{g\}$;
 \EndWhile
 \State\Return{$\id{h}$};
\EndProcedure
\end{algorithmic}
\end{algorithm}

\begin{algorithm}[H]
\caption{\label{alg:ed}Computing the embedding dimension of $\St_h(\id{j})$}
\begin{algorithmic}[1]
 \Procedure{EmbeddingDimension}{$\id{j}$,$\prec$}
 \State Execute lines {\small\textsc{2 -- 8}} of \Call{Gr\"obnerStratum}{$\id{j}$,$\prec$}
 \State $L \leftarrow 0$;
 \ForAll{$(\ldots,X^{\delta_i},\ldots,X^{\delta_j},\ldots) \in \Syz(\id{j})$}
 \State $S \leftarrow X^{\delta_i}F_{\alpha_i} - X^{\delta_j}F_{\alpha_j}$;
 \State Compute the reduction $S_{\textnormal{red}}$ of $S$ w.r.t. the \emph{monomial ideal} $\id{j}$;
 \State Collect the set $\mathcal{R}$ of the \emph{linear} $X$-coefficients of $S_{\textnormal{red}}$;
 \State $L \leftarrow L + \langle\mathcal{R}\rangle$;
 \EndFor
 \State \Return{$(\dim k[C]_1 - \dim L)$};
\EndProcedure
\end{algorithmic}
\end{algorithm}

\subsection*{Example: $\hilb_{4z}^3$} 
In \cite{gotzmann-2008}, Gotzmann consider the complete list of the Borel-fixed, saturated, monomial ideals of $k[X_0,X_1,X_2,X_3]$ corresponding to points of $\hilb_{4z}^{3}$. They are:
\[
\begin{array}{l}
\id{b}_3 = (X_0^2,X_0 X_1,X_1^3),\\
\id{b}_4 = (X_0^2,X_0 X_1,X_0 X_2^2,X_1^4),\\ 
\id{b}_5 = (X_0^2,X_0 X_1,X_0 X_2,X_1^5,X_1^4 X_2), \\
\id{b}_6 = (X_0,X_1^5,X_1^4 X_2^2).
\end{array}
\]
The  index $s$ in $\id{b}_s$ is the regularity of the ideal. Moving from this point, Gotzmann proves that there are two irreducible components: the first containing $\id{b}_3$ with dimension 16 and the second (the Reeves-Stillman one) containing $\id{b}_6$ with dimension 23. Here we obtain a computational confirmation of this result. Furthermore we also prove that the two components are reduced, rational and that they have a transversal intersection.

Since the Gotzmann number of the Hilbert polynomial $p(z) = 4z$ is $6$, to deduce informations about $\hilb_{4z}^3$ using the results obtained in \S\ \ref{sec:opensubset}, we have to consider the truncated ideals $\id{j}_s = (\id{b}_s)_{\geqslant 6}$. 

In the case $\textnormal{char}\, k = 0$, a Borel ideal $\id{i}$ is characterized by the combinatorial property
\begin{equation}\label{borelrel}
 X_i X^\alpha \in \id{i}\quad\Longrightarrow\quad X_{i-1} X^\alpha \in \id{i}.
\end{equation}
As shown in \cite{sherman-2005}, the set of monomials of a fixed degree of a Borel ideal $\id{i}$ is a filter for the transitive closure of the partial ordering $\leq_B$ induced by the relation \eqref{borelrel} (${X^\alpha}{X_{i-1}} >_B X^\alpha X_{i}$). For each $s$, we can look for a term ordering $\prec_s$, obtained refining the partial order $\leq_B$, such that the ideal $\id{j}_s$ becomes a $(6,\prec_s)$-segment. 

It is possible to achieve this result considering a term ordering given by a matrix of the type:
\[
\left(\begin{array}{cccc} 
1&1&1&1\\ 
a&b&c&d\\
0&1&0&0\\ 
0&0&1&0
\end{array}\right).
\]
More precisely:
\begin{itemize}
	\item $\id{j}_3$ is a segment w.r.t. $\prec_3$ given by  $(a,b,c,d)=(3,2,1,1)$;
	\item $\id{j}_4$ is a segment w.r.t. $\prec_4$ given by  $(a,b,c,d)=(15,5,2,1)$;
	\item $\id{j}_5$ is a segment w.r.t. $\prec_5$ given by   $(a,b,c,d)=(9,3,2,1)$;
	\item $\id{j}_6$ is a segment w.r.t. $\prec_6$ given by   $(a,b,c,d)=(1,0,0,0)$;
\end{itemize}

Let us now examine these ideals one at the time.

\noindent $\id{j}_6$. It is the lexsegment ideal $\id{L}(2,4,0)$, that is $\prec_6$ is $\mathtt{DegLex}$ (or $\mathtt{Lex}$, which is the same for homogeneous ideals). In the previous section we proved that the Gr\"obner stratum $\St_h(\id{L}(2,4,0),\texttt{Lex})$ is an open subset of $\hilb_{4z}^3$ and that it is an affine space; an easy computation gives the dimension 23.

\medskip 

\noindent $\id{j}_3$. As $X_2$ does not appear in the monomials of the basis of $\id{b}_3$, thanks to Theorem \ref{th:saturato}, we know that $\St_h(\id{j}_3,\prec_3) \simeq \St_h(\id{b}_3,\prec_3)$. Using a computer procedure based on Algorithm \ref{alg:idSth} we obtain  explicit equations for $\St_h(\id{b}_3,\prec_3)$, finding that it is an affine space of dimension 16, which is isomorphic to an open subset of the component named $H_{\textit{VA}}$ by Gotzmann (after Vainsencher-Avritzer \cite{VA}).

\medskip

\noindent $\id{j}_4$. Always by Theorem \ref{th:saturato}, we can compute $\St_h(\id{b}_4,\prec_4)$ instead of $\St_h(\id{j}_4,\prec_4)$. Using the same computer procedure as below, we compute that initially there are $44$ variables $C$, but $20$ of them can be eliminated and the minimal embedding realizes $\St_h((\id{b}_4),\prec_4)$ as an affine subscheme $S$ of $\mathbb{A}^{24}$ given by an ideal $\id{h}$. Moreover $\id{h}$ is the product of a principal ideal $(K)$ (more precisely $K$ is the coefficient of the monomial $X_1^3$ in the polynomial $F_i$ such that $\Lt(F_i)=X_0X_2^2$) and an ideal $\id{h}_1$. The ideal $(K)$ defines an hyperplane in $\mathbb{A}^{24}$, which is an open subscheme of the $H_{RS}$ component (its dimension is 23). The other ideal $\id{h}_1$ defines an open subscheme of $H_{AV}$.  Looking at $\id{h}_1$ it is possible to see that there are some more eliminable variables and that the minimal embedding gives an isomorphism with  $\mathbb{A}^{16}$. Looking at the ideals $(K)$ and $\id{h}_1$, it is easy to check that they have a transversal intersection that is the hyperplane defined by $K$ does not contain the Zariski tangent space to $H_{AV}$  at each point in $H_{RS}\cap H_{AV}$ . 

\medskip

\noindent $\id{j}_5$. Applying Theorem \ref{th:saturato}, we compute the structure of the open subset $\HH_{\id{b}_5}$ of $\hilb_{4z}^3$ by the computation of $\St_h((\id{b_5})_{\geqslant 4},\prec_5)$. In this case there are 344 new variables $C$: 317 are eliminable, so that the embedding dimension is 27, that is the point corresponding to $\id{b}_5$ is singular in the Hilbert scheme $\hilb_{4z}^3$. Going through the computation, we find a Gr\"obner basis of the ideal $\id{h}((\id{b_5})_{\geqslant 4})$ defined by 9 polynomials. The stratum is irreducible (and so there are no new components): indeed the open subset obtained excluding the hyperplane  defined by the coefficient of the monomial $X_1^4$ of the polynomial with leading term $X_0X_2X_3^2$ is isomorphic to an open subset of an affine space of dimension 23, that is to an open subset of the Reeves-Stillman component. 
 On the other hand cutting the stratum with this same hyperplane we obtain an equi-dimensional  subscheme of dimension 22 and the same degree than $\St_h((\id{b_5})_{\geqslant 4},\prec_5)$, which is scheme theoretically the union of two irreducible components $V_1$ and $V_2$,
 with the same Hilbert polynomial. One of them  can be naturally identified with the  stratum of the saturated ideal $\id{b_5}$ in the sense that   they are isomorphic and moreover  their points correspond to the same curves in $\PP^3$ (see Theorem \ref{th:saturato} \emph{\ref{it:saturato_ii})}):  this component  is obtained cutting $\St_h((\id{b_5})_{\geqslant 4},\prec_5)$ with the hyperplanes defined by the coefficients of $X_1^4$ in the  polynomials with leading term $X_0X_2X_3^2$, $X_0X_1X_3^2$ and $X_0^2X_3^2$ respectively (see Theorem \ref{th:saturato} \emph{\ref{it:saturato_v}) } and \emph{\ref{it:saturato_vv})}).  The other one $V_2$ can be obtained from $V_1$ up to a special change of coordinates in $\PP^3$. Finally we can verify that the singular locus of $\St_h((\id{b_5})_{\geqslant 4},\prec_5)$ is contained in the intersection of $V_1$ and $V_2$.


\bigskip

\noindent Paolo Lella

\noindent {\small paolo.lella@unito.it}

\noindent Dipartimento di Matematica dell'Universit\`{a} di Torino\\ 
         Via Carlo Alberto 10 \\ 
         10123 Torino, Italy  
  
  \medskip
        
   \noindent Margherita Roggero

    \noindent {\small margherita.roggero@unito.it}      

\noindent Dipartimento di Matematica dell'Universit\`{a} di Torino\\ 
         Via Carlo Alberto 10 \\ 
         10123 Torino, Italy 

\end{document}